\documentclass[11pt,draftclsnofoot, onecolumn]{IEEEtran}
\usepackage{amsfonts,amssymb,amsmath}
\usepackage{graphicx}
\usepackage{booktabs}
\usepackage{mathtools}
\usepackage{xcolor,color}
\usepackage{amsthm}
\usepackage{algorithm}
\usepackage{algpseudocode}
\usepackage{bbm}
\usepackage{booktabs}
\usepackage{cite}
\allowdisplaybreaks

\DeclarePairedDelimiter{\floor}{\lfloor}{\rfloor}
\newtheorem{theorem}{Theorem}
\newtheorem{lemma}[theorem]{Lemma}
\newtheorem{assumption}[theorem]{Assumption}
\newtheorem{definition}[theorem]{Definition}

\newtheorem{remark}[theorem]{Remark}

\newtheorem{proposition}[theorem]{Proposition}

\begin{document}
\title{Data-Driven Robust Control Using Prediction Error Bounds Based on Perturbation Analysis}

\author{Baiwei Guo, Yuning Jiang, Colin N. Jones,  Giancarlo Ferrari-Trecate
\thanks{This work has been supported by the Swiss National Science Foundation under the NCCR Automation (grant agreement 51NF40\_80545).}
\thanks{
Baiwei Guo and Giancarlo Ferrari-Trecate are with the DECODE
group, Institute of Mechanical Engineering, EPFL, Switzerland.  email:\{\tt baiwei.guo, giancarlo.ferraritrecate@epfl.ch\}}
\thanks{Yuning Jiang and Colin N. Jones are with the PREDICT
group, Institute of Mechanical Engineering, EPFL, Switzerland.  email:\{\tt yuning.jiang, colin.jones@epfl.ch\}
}
}

\maketitle 
\begin{abstract}
For linear systems, many data-driven control methods rely on the behavioral framework, using historical data of the system to predict the future trajectories. However, measurement noise introduces errors in predictions. When the noise is bounded, we propose a method for designing historical experiments that enable the computation of an upper bound on the prediction error. This approach allows us to formulate a minmax control problem where robust constraint satisfaction is enforced. We derive an upper bound on the suboptimality gap of the resulting control input sequence compared to optimal control utilizing accurate measurements. As demonstrated in numerical experiments, the solution derived by our method can achieve constraint satisfaction and a small suboptimality gap despite the measurement noise.
\end{abstract}

\section{Introduction}
With advancements in sensing, storage, and computation technologies, the availability of data has significantly increased. This enables the design of controllers that are both safe (i.e., satisfying constraints) and near-optimal in the presence of uncertainties \cite{tsiamis2022statistical,karimi2017data}.

Control synthesis methods can be broadly categorized as either \textit{model-based} or \textit{model-free}. Model-based methods involve two stages: collecting historical system trajectories to identify a parametric model for system dynamics, and utilizing this model to calculate optimal control decisions. Classical methods like LQG control \cite{athans1971role} and feedback linearization \cite{khalil2015nonlinear} fall under this category. A common drawback is that the involved system identification tasks can be demanding especially when the model structure is unknown, requiring time-consuming procedures such as data preprocessing, model selection, and validation \cite{ljung1999system}. Additionally, regarding the accuracy of the derived model under measurement noise, existing upper bounds for identification errors rely on the unknown ground-truth system model instead of empirical data \cite{dean2020sample,oymak2019non}. Therefore, without further assumptions, these bounds cannot be directly used in robust control.

In contrast, model-free methods derive the optimal control inputs directly from data without explicitly identifying the system model, simplifying the formulation process.  For example, the authors of \cite{yin2021maximum} use a maximum likelihood framework to obtain an optimal data-driven
model for trajectory prediction. This way, one can avoid separate ML-based model identification and Kalman Filter for state estimation. Model-free methods are also highly effective in robust constraint satisfaction~\cite{lian2021adaptive} and optimal control of some nonlinear dynamics \cite{coulson2019data}. However, a complete theoretical explanation of these advantages are still lacking.

Many model-free methods are based on Willems' fundamental lemma \cite{willems2005note} for linear systems. This lemma states that, under the assumption of persistent excitation, any trajectory of a given linear system is a linear combination of a Hankel matrix's columns where the Hankel matrix consists of the historical trajectories. It allows us to represent the system directly with the data instead of identifying a parametric model. The system representation based on the Hankel matrix is called a behavioral model.

In this paper, we develop a behavioral-model-based approach toward robust control of linear systems. Specifically, we aim to minimize the worst-case regulation cost and satisfy constraints under bounded output measurement noise. Several existing works approximately solve this problem from different perspectives and under different assumptions. For input-output systems whose states are not directly measurable, \cite{furieri2022near} optimizes an output feedback controller by using the behavioral model to identify the system's impulse response. The derived controller results in a higher cost than the optimal one due to the measurement noise and the authors derive a bound for the suboptimality gap. The drawback of this method is that both the suboptimality bound and the satisfaction of the constraints rely on a bootstrap procedure which requires extra resampling and computation. Moreover, a rigorous analysis on the accuracy of bootstrap results using finite samples is still lacking \cite{kilian1999finite}. In contrast, without the efforts of identification, \cite{xu2021data} minimizes the worst-case trajectory tracking cost by reformulating the minmax problem into a Semidefinite Program through the S-lemma. However, the optimality results in \cite{xu2021data} rely on the assumption that the noise sequence satisfies a cumulative quadratic constraint. Therefore, they do not hold when all entries of the noise sequence are known to satisfy box constraints. To address this issue, \cite{huang2021robust} adopts a different reformulation where model misfit penalty is added to the tracking cost. The involved minmax problem is solvable by robust optimization techniques and the authors derive a suboptimality gap bound. However, this bound is conservative in the sense that it does not vanish when the noise decreases to zero. 

To develop a method that enjoys both a less conservative suboptimality gap and the guarantee of constraint satisfaction, we leverage perturbation analysis for assessing the influence of measurement noise on the behavioral-model predictions.
Research on perturbation analysis of data-driven prediction includes \cite{berberich2020data}, \cite{coulson2022robust} and \cite{coulson2022quantitative}. In \cite{berberich2020data}, the authors use a cost function similar with the one in \cite{huang2021robust}, aiming to minimize the sum of the misfit penalty and the tracking cost. Therefore, the prediction error upper bound relies on a variable related to the optimal control task and thus can only be evaluated after the optimal control problem has been solved. In contrast, the robust fundamental lemma proposed in \cite{coulson2022robust} gives a prediction error analysis independent of control tasks, but it is limited to input-state systems and one-step-ahead prediction. The work in \cite{coulson2022quantitative} extends \cite{coulson2022robust} by deriving lower bounds for the singular values of the Hankel matrix in the behavioral model for input-output systems. However, the lower bounds depend on the unknown ground truth system model.

In this paper, we aim to upperbound the error of behavioral-model-based prediction  for general input-output linear systems and utilize this bound for minmax robust control with guarantees of constraint satisfaction. Our main contributions are the following:
\begin{itemize}
\item We propose an input generation strategy to collect historical data under bounded measurement noise for the construction of a Page matrix (a variant of Hankel matrix, see \cite{coulson2020distributionally}) enabling the derivation of an error bound of the data-driven predictions, which only relies on the noisy data. This bound is valid when the historical inputs achieve sufficient ``persistent-excitation-to-noise ratio'' (rigorously stated in Assumption \ref{ass: small_noise}) and the observability index (see Definition \ref{def: observability index}) is identified correctly. The first condition can be satisfied if collecting multiple historical data sets for averaging the noisy measurements or enlarging the input signals is allowed. We achieve the second through a data-driven method with correctness guarantee.
\item For unconstrained regulation of MISO systems, in order to minimize the worst-case cost, we utilize the new prediction error bound to formulate a minmax problem and bound the suboptimality gap. The derived bound decreases to zero as the measurement noise converges to zero. This scheme can be extended to regulation of MIMO systems and robust constraint satisfaction.

\end{itemize}

The rest of this paper is organized as follows: in Section \ref{sec: prediction}, we introduce the basics of data-driven control and formulate the robust control problem. Section \ref{sec: datadrivenprediction} focuses on MISO systems with no input output constraints and addresses perturbation analysis on the noisy behavioral model used for trajectory prediction. The robust control scheme is proposed in Section \ref{sec: formulationofGu} where the associated upper bound of the worst-case cost derived by using perturbation analysis is minimized and the suboptimality gap is bounded. We extend our method for regulation of MIMO systems and for robust constraint satisfaction in Section V. Experiments illustrating the performance of the proposed approach are shown in Section~\ref{sec: numerical studies}.

\textit{Notations:} Given a time-varying vector variable $v$, we use $v_t$ to denote its value at time instant $t$, let $[t_1,t_2] = \{t_1,t_1+1,...,t_2\}$ and set $v_{[t_1,t_2]}: =\{v_{t_1},v_{t_1+1},\ldots,v_{t_2}\}$, 
$\mathrm{col}(v_{[t_1,t_2]}): = \begin{bmatrix} v^\top_{t_1} & v^\top_{t_1+1} & \ldots & v^\top_{t_2}\end{bmatrix}^\top$. We use $\|x\|_i$ to denote the $\ell_i$ norm of $x$. Moreover, $||x||$ is the $\ell_2$ norm. Given positive semidefinite matrix $Q$, the term $\|x\|^2_Q$ denotes $x^\top Q x.$ For a matrix $M$, $||M||_i$ denotes the matrix $i$-norm while $||M||_\mathrm{max}: = \max_{i,j}|m_{ij}|$ and $||M||:= ||M||_2$. We also denote $\sigma_{\min}(M)$ as the smallest singular value of $M$ and $\sigma_{\max}(M)$ as the largest. We use $M_{i,\cdot}$ to denote the $i$-th row, $M_{\cdot,i}$ the $i$-th column and $M_{i:j,\cdot}$ the submatrix consisting of the rows of $M$ from the $i$-th to the $j$-th. For a given $x\in\mathbb{R}$, we use notation $\floor{x}$ to denote the floor function, i.e., $\floor{x} = \max\{z\in\mathbb{Z}\mid z\leq x\}$. For $y\in\mathbb{R}^n$ and $r>0$, we let $\mathcal{B}(y,r):=\{x:\|x-y\|\leq r\}$. The identity matrix with $n$ rows is denoted as $I_n$. The pseudo-inverse of a matrix $H$ is written as $H^\dagger$.

\section{Preliminaries and Problem Formulation}
\label{sec: prediction}
\subsection{Preliminaries: data-driven description of linear systems}
This paper considers the regulation of a discrete-time linear time-invariant (LTI) system with the following controllable and observable minimal realization,
\begin{equation}\label{eq:LTI_system_ss} 
x_{t+1} =A x_{t}+B u_{t}\,,\;\;
y_{t} =C x_{t}+D u_{t}, \tag{$*$}
\end{equation}
where $u_t\in \mathbb{R}^m$, $y_t\in \mathbb{R}^p$ and  $x_t\in \mathbb{R}^{n_x}$. We assume that $n_x$ and the system matrices $A,B,C,D$ are unknown. Instead of identifying the system matrices, we utilize a behavioral model to characterize the possible trajectories in a horizon of length $L$. To build this system representation, we excite the linear system with an input sequence $u_{[1,T]}$ and collect the output data $y_{[1,T]}$, where $T>L$. The
$L\text{-Page}$ matrix of $u_{[1,T]}$ is given by
\[
\mathcal{P}_L(u_{[1,T]}): = \begin{bmatrix}
u_1 & u_{L+1} & \cdots & u_{l_\mathrm{h}L-L+1} \\
u_2 & u_{L+2} & \cdots & u_{l_\mathrm{h}L-L+2} \\
\vdots & \vdots & \ddots & \vdots \\
u_L & u_{2L} & \cdots & u_{l_\mathrm{h}L},
\end{bmatrix}
\]
where $l_\mathrm{h} = \floor{\frac{T}{L}}$ denotes the number of columns \cite{coulson2020distributionally}.
To evaluate whether the input sequence $u_{[1,T]}$ along with the resulting outputs is sufficiently informative to uniquely determine system \eqref{eq:LTI_system_ss}, we introduce the following definition analogous to persistent excitation in system identification.
\begin{definition}[\cite{coulson2020distributionally}]
\label{def: persistent_excitation}
For $L, T, d\in \mathbb{Z}^+$, we say the input sequence $u_{[1,T]}$ is $L$-Page exciting of order $d$ if the following matrix has full row rank,  
$$
\mathcal{P}_{L,d}\left(u_{[1,T]}\right):=\left(\begin{array}{c}
\mathcal{P}_{L}\left(u_{[1, T-(d-1) L]}\right) \\
\mathcal{P}_{L}\left(u_{[L+1, T-(d-2) L]}\right) \\
\vdots \\
\mathcal{P}_{L}\left(u_{[L(d-1)+1, T]}\right)
\end{array}\right).
$$
\end{definition}
By using the collected input output data $(u_{[1,T]},y_{[1,T]})$, called \textit{historical}, one might be able to determine whether another trajectory $(u^{r}_{[1,L]},y^{r}_{[1,L]})$, called \textit{recent}, is generated by system \eqref{eq:LTI_system_ss}. Rigorously, we have the following result, which is a variant of the well-known Willem's Fundamental Lemma \cite{willems2005note}.

\begin{lemma}[{\cite[Theorem~2.1]{coulson2020distributionally}}]
\label{lmm: Willem's lemma}
For the LTI system described in \eqref{eq:LTI_system_ss}, given a $T$-length historical trajectory $(u_{[1,T]},y_{[1,T]})$ where $u_{[1,T]}$ is $L$-Page exciting of order $n_{x}+1$ and the $L$-length recent trajectory $(u^{r}_{[1,L]},y^{r}_{[1,L]})$, there exists $x^{r}_{[1,L]}$ with $x_i\in \mathbb{R}^{n_x}$, $1\leq i \leq L$, such that $(x^{r}_{[1,L]},u^{r}_{[1,L]},y^{r}_{[1,L]})$ satisfies \eqref{eq:LTI_system_ss} if and only if there exists a vector $g\in \mathbb{R}^{l_\mathrm{h}}$ such that 
\begin{equation}
\label{eq: willem_lemma}
\begin{bmatrix}
\mathcal{P}_{L}\left({u}_{[1, T]}\right) \\
\mathcal{P}_{L}\left({y}_{[1, T]}\right)
\end{bmatrix}
g=\begin{bmatrix}
\mathrm{col}(u^{r}_{[1, L]}) \\
\mathrm{col}(y^{r}_{[1, L]})
\end{bmatrix}.
\end{equation}
\end{lemma}

\subsection{Problem formulation: robust control under measurement noise}
Given the historical data $(u_{[1,T]},y_{[1,T]})$ where $u_{[1,T]}$ is $L$-Page exciting of order $n_x+1$ and an $l_\mathrm{p}$-long initial trajectory $(u^{r}_{[1,l_\mathrm{p}]},y^{r}_{[1,l_\mathrm{p}]})$ with $l_\mathrm{p}<L$, we consider the following regulation problem for the trajectory $(u^{r}_{[1,T]},y^{r}_{[1,T]})$ from $t = l_\mathrm{p}+1$ to $t =L$:
\begin{equation}
\label{eq: original_problem}
\begin{aligned}
\min_{u^{r}_{[l_\mathrm{p}+1,L]},y^{r}_{[l_\mathrm{p}+1,L]}} & \sum_{i=1}^{L-l_\mathrm{p}}(|| u^{r}_{l_\mathrm{p}+i}||^2 + ||y^{r}_{l_\mathrm{p}+i}||^2)\\
\quad\mathrm{s.t.}  \quad \quad & \text{there exists $g$ such that} \\
& u_{[1,T]},y_{[1,T]},u^{r}_{[1,L]},y^{r}_{[1,L]} \text{ satisfy \eqref{eq: willem_lemma}}.
\end{aligned}
\end{equation}
For convenience, we let $l_\mathrm{f} = L-l_\mathrm{p}$ and use the following notations for historical data,
\begin{equation}
\label{eq: Page_matrices}
\begin{aligned}
U_\mathrm{p} & = \begin{bmatrix}
I_{ml_\mathrm{p}} & 0
\end{bmatrix} \mathcal{P}_{L}\left({u}_{[1, T]}\right), U_\mathrm{f} = \begin{bmatrix}
0 & I_{ml_\mathrm{f}}
\end{bmatrix} \mathcal{P}_{L}\left({u}_{[1, T]}\right),\\
Y_\mathrm{p} & = \begin{bmatrix}
I_{pl_\mathrm{p}} & 0
\end{bmatrix} \mathcal{P}_{L}\left({y}_{[1, T]}\right), Y_\mathrm{f} = \begin{bmatrix}
0 & I_{pl_\mathrm{f}}
\end{bmatrix} \mathcal{P}_{L}\left({y}_{[1, T]}\right)
\end{aligned}
\end{equation}
and for recent data,
\begin{equation}
\label{eq: split_of_recent}
\begin{aligned}
u_\mathrm{p} & = \mathrm{col}(u^{r}_{[1,l_\mathrm{p}]}),  u_\mathrm{f}  = \mathrm{col}(u^{r}_{[l_\mathrm{p}+1,L]}),\\
y_\mathrm{p} & = \mathrm{col}(y^{r}_{[1,l_\mathrm{p}]}),  y_\mathrm{f}  = \mathrm{col}(y^{r}_{[l_\mathrm{p}+1,L]}).
\end{aligned}
\end{equation}
In practice, output measurements are subject to noise. Specifically, for any $i\in[1,\ldots,T]$, $j\in[1,\ldots,l_\mathrm{p}]$, there exist noise vectors $w_i, w^r_j\in\mathbb{R}^{p}$ such that the measurements are $\hat{y}_i = {y}_i+w_i$ and  $\hat{y}^r_j = {y}^r_j+w^r_j$. We build $\widehat{Y}_\mathrm{p}, \widehat{Y}_\mathrm{f}$ and $ \hat{y}_\mathrm{p}$  from $\hat{y}_{[1,T]}$ and $\hat{y}^r_{[1,l_\mathrm{p}]}$ as noisy counterparts of ${Y}_\mathrm{p}, {Y}_\mathrm{f}$ and $ {y}_\mathrm{p}$. In this paper, we only consider bounded noise, as stated in the following assumption.
\begin{assumption}
\label{ass: measurement_noises}
For the noisy measurements $\widehat{Y}_\mathrm{p}, \widehat{Y}_\mathrm{f}, 
\hat{y}_\mathrm{p},$ corresponding to ${Y}_\mathrm{p}, {Y}_\mathrm{f}, 
{y}_\mathrm{p},$ the following holds,
$$\mathrm{max}\{||\widehat{Y}_\mathrm{p} - {Y}_\mathrm{p}||_\mathrm{max},||\widehat{Y}_\mathrm{f} - {Y}_\mathrm{f}||_\mathrm{max},||\hat{y}_\mathrm{p} - {y}_\mathrm{p}||_\infty\}\leq \delta,$$
where $\delta>0$ is a known constant.
\end{assumption}
Given a fixed control strategy, different noise realizations lead to different performances. To attenuate the influence of the uncertainties, we aim to \textbf{design a robust control scheme where the worst-case regulation cost is minimized}. In general, the worst-case cost is hard to compute and we will provide an upper bound $c_{\textrm{worst}}({u}_\mathrm{f})$ to it given the historical data $u_{[1,T]},\hat{y}_{[1,T]}$ and the recent data $u_\mathrm{p}$ and $\hat{y}_\mathrm{p}$. We formulate the robust control problem as $$u^*_\mathrm{f}=\mathrm{argmin}_{u_\mathrm{f}} c_{\textrm{worst}}(u_\mathrm{f}).$$ This way, we ensure that the true cost induced by $u^*_\mathrm{f}$ is less than $c_{\textrm{worst}}(u^*_\mathrm{f})$. 

In Section \ref{sec: formulationofGu} we propose a formulation of $c_{\textrm{worst}}(u_\mathrm{f})$. But before that, in Section \ref{sec: datadrivenprediction}, by assuming $u_\mathrm{f}$ is given, we introduce tools of data-driven prediction based on behavioral models and conduct perturbation analysis, which allows us to estimate the worst-case cost by using the noisy data. To avoid bulky statements, we {only} consider the MISO case in Sections \ref{sec: datadrivenprediction} and \ref{sec: formulationofGu}, while the extension to the MIMO case is presented in Section \ref{sec: computation+extensions}.

\section{Data-Driven Prediction with Perturbation Analysis for MISO Systems}
\label{sec: datadrivenprediction}
In this section, we investigate, for MISO systems, a data-driven prediction method that generates an estimate $\hat{y}_\mathrm{f}$ for ${y}_\mathrm{f}$ based on historical data along with $u_\mathrm{p}$, $\hat{y}_\mathrm{p}$ and $u_\mathrm{f}$. Specifically, we propose a method for generating historical data enabling the derivation of an upper bound for the prediction error. This ingredient is essential for the computation of upper bounds to worst-case costs in Section \ref{sec: formulationofGu}.

\subsection{The proposed data-driven prediction scheme}
For MISO systems, we adopt the data-driven prediction method used in the framework of \textbf{PEM-MPC}\footnote{PEM stands for \textit{Prediction Error Method}.} proposed in~\cite{huang2019data}. Specifically, we look into the following prediction scheme
\begin{equation}
\label{eq: ls_pred_noisy}
\hat{g}(u_\mathrm{f}) = \widehat{H}^\dagger 
\hat{b}(u_\mathrm{f}),\;
\hat{y}_\mathrm{f}(u_\mathrm{f})  = \widehat{Y}_\mathrm{f}\hat{g}(u_\mathrm{f}),\text{ where }\widehat{H} = \begin{bmatrix}
U^\top_\mathrm{p} & \widehat{Y}^\top_\mathrm{p} & U^\top_\mathrm{f}
\end{bmatrix}^\top, \;\hat{b}(u_\mathrm{f})= \begin{bmatrix}
u^\top_\mathrm{p} & \hat{y}^\top_\mathrm{p} & u^\top_\mathrm{f}
\end{bmatrix}^\top.
\end{equation}
We also define the noiseless counterparts of $\widehat{H}, \hat{b}(u_\mathrm{f}),\hat{g}(u_\mathrm{f}),\hat{y}_\mathrm{f}(u_\mathrm{f})$ in \eqref{eq: ls_pred_noisy} as $\overline{H}, \bar{b}(u_\mathrm{f}),\bar{g}(u_\mathrm{f}),\bar{y}_\mathrm{f}(u_\mathrm{f})$. One can easily verify that $\|\bar{g}(u_\mathrm{f})\| = \min\{\|g\|:\bar{H}{g} = \bar{b}(u_\mathrm{f}) \text{ and } {Y}_\mathrm{f}g = \bar{y}_\mathrm{f}(u_\mathrm{f})\}$. Before elaborating on the historical data collected and the construction of $\widehat{H}$ in Assumption \ref{ass: design_of_historical_data}, we introduce the concept of observability index.
\begin{definition}
\label{def: observability index}
We say $l_\mathrm{o}$ is the observability index of the system \eqref{eq:LTI_system_ss} if the $l_\mathrm{o}$ is the smallest positive integer $l$ such that the observability matrix
$\mathcal{O}(l) := 
\begin{bmatrix}
C^\top &
(C A)^\top &
\cdots &
(C A^{l-1})^\top
\end{bmatrix}^\top
$ is full column rank.
\end{definition}
\begin{assumption}
\label{ass: design_of_historical_data}
Given horizon length $L>l_o$, the historical inputs of the MISO system \eqref{eq:LTI_system_ss} are generated using the following setting:
\begin{subequations}
\begin{align}x_{1}=0,\\
u_{[1,L]},u_{[L+1,2L]},\ldots,u_{[4L^2+1,4L^2+L]} \overset{\mathrm{i.i.d.}}{\sim}\mathcal{Q},\\
T\geq 4L^2+L \text{ and } u_{[1,T]} \text{ is $L$-Page exciting of order $n_x+1$,}
\label{eq:iid_input}
\end{align}
\label{eq: design_of_inputs}
\end{subequations}\vspace{-0.8cm}\\
where $\mathcal{Q}$ is a probability distribution such that for $x\in \mathbb{R}^{mL}$ 
\begin{equation}
\label{eq: zero_mass}
\text{if } x\sim \mathcal{Q} \text{ then for any subspace }V\subsetneqq \mathbb{R}^{mL},\text{ } \mathbb{P}(\{x\in V\})=0.   
\end{equation}
The Page matrices $U=\mathcal{P}_{L}\left({u}_{[1, T]}\right)$ and $\widehat{Y}=\mathcal{P}_{L}\left(\hat{y}_{[1, T]}\right)$ are split into $U_\mathrm{p}, U_\mathrm{f},\widehat{Y}_\mathrm{p}, \widehat{Y}_\mathrm{f}$ using \eqref{eq: Page_matrices} with 
\begin{equation}
l_\mathrm{p} = l_\mathrm{o},
\end{equation}
and $\widehat{H}$ is built using \eqref{eq: ls_pred_noisy}.
\end{assumption}
For identification of $l_\mathrm{o}$ using data, we refer to Section \ref{sec: identification_l_o} and Algorithm~\ref{alg:observability index identification}. Under Assumption \ref{ass: design_of_historical_data}, we will see in Section \ref{sec: perturbation_analysis_ls} that the noiseless matrix $\overline{H}$ is almost surely full row rank. This property is essential for upperbounding the prediction error resulting from \eqref{eq: ls_pred_noisy}. The reason is that, if $\overline{H}$ has a zero singular value, even infinitesimal measurement noise can result in a huge prediction error due to the use of the pseudoinverse in \eqref{eq: ls_pred_noisy}.

\subsection{Perturbation analysis}
\label{sec: perturbation_analysis_ls}
We show in Theorem \ref{thm: full-rankness} and Theorem~\ref{thm: prediction_error_bound} that the prediction error resulting from the scheme \eqref{eq: ls_pred_noisy} can be bounded. 
\begin{theorem}
\label{thm: full-rankness}
Under Assumption \ref{ass: design_of_historical_data}, we have $\mathbb{P}(\{\overline{H} \text{ is full row rank}\}) = 1$.
\end{theorem}
The proof of Theorem \ref{thm: full-rankness} is given in Appendix \ref{sec: proof_full_rank}.

\begin{remark}
\label{rmk: mimonotgood}
Theorem \ref{thm: full-rankness} does not hold for MIMO systems. As a simple example, when there are two outputs measuring the same quantity, almost surely the noiseless $\overline{H}$ is not full row rank.  
\end{remark}

By utilizing Theorem \ref{thm: full-rankness}, we show in Theorem \ref{thm: prediction_error_bound} that, when the measurement noise is sufficiently small (as indicated below), the prediction errors given by the least-square solution \eqref{eq: ls_pred_noisy} can be bounded.
\begin{assumption} \label{ass: small_noise}
With $\widehat{H}$ constructed according to Assumption \ref{ass: design_of_historical_data}, we have
\begin{equation}
\label{eq: small_noise}
    \delta<\frac{\sigma_{\min}(\widehat{H})}{2l_{\mathrm{h}}}.
\end{equation}
\end{assumption}

\begin{remark}
\label{rmk: signal_to_noise_ratio}
As will be shown later in Theorem \ref{thm: prediction_error_bound}, under Assumption \ref{ass: small_noise}, an upper bound for the prediction error is roughly proportional to $\sigma_{\min}(\widehat{H})^{-1}$. Therefore, we can regard $\sigma_{\min}(\widehat{H})$ as a measure of ``persistent excitation'' and thus \eqref{eq: small_noise} requires the ``persistent-excitation-to-noise ratio'' to be sufficiently large. In \cite[Lemma 1]{berberich2020data}, a similar assumption is made. To satisfy Assumption \ref{ass: small_noise}, one can decrease the noise magnitude by collecting multiple historical data sets and using averaging techniques if the entries of the noise sequence are independent and identically distributed (for details see Remark \ref{rmk: averaging}). Since from \cite[Theorem 4.3]{li1998relative} we have
\begin{equation}
\label{eq: singular_value_perturbation}
\sigma_{\min}(\widehat{H})\geq \sigma_{\min}(\overline{H})-\|\widehat{H}-\overline{H}\|,   
\end{equation} we can also amplify the historical input signals in \eqref{eq: design_of_inputs} for obtaining a larger $\sigma_{\min}(\overline{H})$ and hence increasing $\sigma_{\min}(\widehat{H})$. 
\end{remark}
\begin{theorem}
\label{thm: prediction_error_bound}
We denote $\overline{\Delta}_g(u_\mathrm{f}):=\overline{g}(u_\mathrm{f})-\hat{g}(u_\mathrm{f})$ and $\overline{\Delta}_{y_\mathrm{f}}(u_\mathrm{f}):=\bar{y}_\mathrm{f}(u_\mathrm{f})-\hat{y}_\mathrm{f}(u_\mathrm{f})$. If Assumption \ref{ass: small_noise} holds, we have for MISO systems $$||\overline{\Delta}_{g}(u_\mathrm{f})||\leq\mathcal{C}(u_{\mathrm{f}})\delta, \text{ where } \mathcal{C}(u_\mathrm{f})=2\sigma_{\min}(\widehat{H})^{-1}(\sqrt{l_{\mathrm{p}}}+l_{\mathrm{h}}||\hat{g}(u_{\mathrm{f}})||), \text{ and}$$
$$||\overline{\Delta}_{y_{\mathrm{f}}}(u_\mathrm{f})||\leq \mathcal{C}(u_\mathrm{f})||\widehat{Y}_{\mathrm{f}}||\delta + l_\mathrm{h} (\|\hat{g}(u_\mathrm{f})\|+\mathcal{C}(u_\mathrm{f}))\delta.$$
\end{theorem}
\begin{proof}
We let $E:=\widehat{H}-\overline{H}$. Due to Assumption \ref{ass: measurement_noises}, we have 
\begin{equation}
\label{eq: upperbound_of_E}
    ||E||\leq l_{\mathrm{h}}\delta.
\end{equation}
Then the following inequalities hold,
\begin{equation}
\label{eq: inequality_on_g}
\begin{aligned}\frac{\sigma_{\min}(\widehat{H})}{2}||\hat{g}({u}_{\mathrm{f}})-\bar{g}({u}_{\mathrm{f}})|| & \leq(\sigma_{\min}(\widehat{H})-||E||)||\hat{g}({u}_{\mathrm{f}})-\bar{g}({u}_{\mathrm{f}})||\\
& \leq\sigma_{\min}(\widehat{H}-E)||\hat{g}({u}_{\mathrm{f}})-\bar{g}({u}_{\mathrm{f}})||\\
& \leq||\overline{H}(\hat{g}({u}_{\mathrm{f}})-\bar{g}({u}_{\mathrm{f}}))||\\
& \leq ||\widehat{H}\hat{g}({u}_{\mathrm{f}})- \overline{H}\bar{g}({u}_{\mathrm{f}}) + (\overline{H} - \widehat{H})\hat{g}({u}_{\mathrm{f}})||
\\
& \leq||\hat{b}({u}_{\mathrm{f}})-\bar{b}({u}_{\mathrm{f}})||+||{E}||\cdot||\hat{g}({u}_{\mathrm{f}})||\\
& ={(\sqrt{l_{\mathrm{p}}}+l_{\mathrm{h}}||\hat{g}({u}_{\mathrm{f}})||)\delta},
\end{aligned}
%\label{eq:Hbound}
\end{equation}
where the first inequality is due to \eqref{eq: small_noise} in Assumption \ref{ass: small_noise}, the second results from \eqref{eq: singular_value_perturbation} and $\hat{b}(u_\mathrm{f})$ is defined in \eqref{eq: ls_pred_noisy}. 
By simplifying \eqref{eq: inequality_on_g}, we have $$||\overline{\Delta}_{g}(u_\mathrm{f})||\leq\mathcal{C}(u_{\mathrm{f}})\delta.$$ Based on this inequality, we have 
\begin{equation}
\begin{aligned}||\hat{y}_{\mathrm{f}}(u_{\mathrm{f}})-\bar{y}_{\mathrm{f}}(u_{\mathrm{f}})|| & =||\widehat{Y}_{\mathrm{f}}(\hat{g}(u_{\mathrm{f}})-\bar{g}(u_{\mathrm{f}}))+(\widehat{Y}_{\mathrm{f}} - {Y}_{\mathrm{f}})\bar{g}(u_{\mathrm{f}})||\\
& \leq||\widehat{Y}_{\mathrm{f}}||\cdot||\hat{g}(u_{\mathrm{f}})-\bar{g}(u_{\mathrm{f}})||+||\widehat{Y}_{\mathrm{f}} - {Y}_{\mathrm{f}}||\cdot||\bar{g}(u_{\mathrm{f}})||\\
& \leq||\widehat{Y}_{\mathrm{f}}||\cdot||\hat{g}(u_{\mathrm{f}})-\bar{g}(u_{\mathrm{f}})||+\\
& \quad\quad l_{\mathrm{h}}\delta(||\hat{g}(u_{\mathrm{f}})||+||\bar{g}(u_{\mathrm{f}})-\hat{g}(u_{\mathrm{f}})||)\\
& \leq \mathcal{C}(u_\mathrm{f})||\widehat{Y}_{\mathrm{f}}||\delta + l_\mathrm{h} (\|\hat{g}(u_\mathrm{f})\|+\mathcal{C}(u_\mathrm{f}))\delta.
\end{aligned}
\label{eq:yfbound_any_u}
\end{equation}
\end{proof}

In the literature, similar results on perturbation analysis accounting for measurement noise can be found in \cite{oymak2019non} and \cite{dean2020sample}. However, the associated bounds utilize the unknown true system model. In the behavioral model framework, \cite{berberich2020data} also derives a prediction error bound after solving an optimal control problem. Therefore, the bound is only valid for the  optimal control sequence $u^*_\mathrm{f}$ determined by the specific problem formulation. In contrast, our upper bound in Theorem \ref{thm: prediction_error_bound} can be calculated directly from the noisy data for any given $u_\mathrm{f}$. This feature allows us to formulate in Sections \ref{sec: formulationofGu} and \ref{sec: computation+extensions} a min-max regulation problem where robust constraint satisfaction can be enforced. 

\subsection{Data-driven identification of $l_\mathrm{o}$}
\label{sec: identification_l_o}
We discuss the identification of $l_\mathrm{o}$, which is needed in the construction of $\widehat{H}$ for satisfying Assumption~\ref{ass: design_of_historical_data}. For this aim, we discuss in Proposition \ref{prop: not-full-rank} whether $\overline{H}$ almost surely has full rank when Assumption \ref{ass: design_of_historical_data} does not hold. Based on this result, we propose a method to identify $l_\mathrm{o}$.

\begin{proposition}
\label{prop: not-full-rank}
If Assumption \ref{ass: design_of_historical_data} does not hold because $l_\mathrm{p}< l_\mathrm{o}$, then$$\mathbb{P}(\{\overline{H} \text{ is full row rank}\}) = 1.$$ Furthermore, if $l_\mathrm{p}> l_\mathrm{o}$, $\mathbb{P}(\{\overline{H} \text{ is full row rank}\}) = 0.$
\end{proposition}
The proof is reported in Appendix \ref{app: proof_of_not_full_row_rank}. Theorem \ref{thm: full-rankness} and Proposition \ref{prop: not-full-rank} say that $l_\mathrm{o}$ is the largest value for $l_\mathrm{p}$ such that the constructed $\overline{H}$ is full row rank. We notice that the noise matrix $E=\widehat{H}-\overline{H}$ satisfies that $\|E\|\leq l_\mathrm{h}\delta$ \footnote{To show this result, one only needs Cauchy–Schwarz inequality and Assumption \ref{ass: measurement_noises}.}.  Therefore, if $\overline{H}$ is not full row rank, i.e., $\sigma_{\min}(\overline{H})=0$, we have
\begin{equation}
\label{eq: singularvalue inequalities}
\sigma_{\min}(\widehat{H})\leq \sigma_{\min}(\overline{H})+\sigma_{\max}(E)\leq l_\mathrm{h}\delta.
\end{equation}
By utilizing \eqref{eq: singularvalue inequalities}, we propose Algorithm \ref{alg:observability index identification} for identification of $l_\mathrm{o}$ and show the correctness of the derived result in Proposition \ref{prop: l_o_identification}.
\begin{proposition}
\label{prop: l_o_identification}
Under Assumption \ref{ass: small_noise}, if $L>l_\mathrm{o}$, Algorithm \ref{alg:observability index identification} returns $l_\mathrm{o}$, the true observability index, almost surely.
\end{proposition}

\begin{proof}
From Assumption \ref{ass: small_noise} and \eqref{eq: singularvalue inequalities}, we know that $\sigma_{\min}(\widehat{H})\geq 2l_\mathrm{h}\delta$ if $l_\mathrm{p}=l_\mathrm{o}$ and $\sigma_{\min}(\widehat{H})\leq l_\mathrm{h}\delta$ if $l_\mathrm{p}>l_\mathrm{o}$. Since Algorithm \ref{alg:observability index identification} terminates when $\sigma_{\min}(\widehat{H})\leq l_\mathrm{h}\delta$ is verified, we only need to show that $\sigma_{\min}(\widehat{H})>l_\mathrm{h}\delta$ almost surely if $l_\mathrm{p}<l_\mathrm{o}$. 

In the remainder of this proof, we denote with $[U_{ \mathrm{p},k}^\top \;U_{ \mathrm{f},k}^\top]^\top$ and $[\widehat{Y}_{ \mathrm{p},k}^\top \;\widehat{Y}_{ \mathrm{f},k}^\top]^\top$ the partitions of the Page matrices $U$ and $Y$, respectively, where $U_{ \mathrm{p},k}$ has $km$ rows and $Y_{ \mathrm{p},k}$ has $kp$ rows. Correspondingly, we write $\widehat{H}_{k}:=[U_{ \mathrm{p},k}^\top\; \widehat{Y}_{ \mathrm{p},k}^\top\; U_{ \mathrm{f},k}^\top]^\top$. We notice that $\widehat{H}_{k}$ is a submatrix consisting of a fraction of $\widehat{H}_{l_\mathrm{o}}$'s rows and $\widehat{H}_{l_o}$ is full row rank almost surely (Theorem \ref{thm: full-rankness}). According to Lemma \ref{lmm: supp_l_o_id} in Appendix \ref{sec: supp_l_o_id}, we have that almost surely $\sigma_{\min}(\widehat{H}_{k})\geq \sigma_{\min}(\widehat{H}_{l_\mathrm{o}})\geq 2l_\mathrm{h}\delta$ for any $k<l_\mathrm{o}$.
\end{proof}
\begin{remark}
\label{rmk: supplement_to_algorithm}
We have another heuristic method which can be used as supplement to Algorithm \ref{alg:observability index identification} for the case when Assumption \ref{ass: small_noise} is not satisfied. Specifically, in Line 1 of Algorighm \ref{alg:observability index identification}, after deriving $(u_{[1,T]},\hat{y}_{[1,T]})$, we can generate several extra historical trajectories $( u(\alpha)_{[1,T]},y(\alpha)_{[1,T]})$ using the initial state $x_1=0$ and the input sequence $u(\alpha)_{[1,T]} = \alpha u_{[1,T]}$ for different values of $\alpha$. Based on the data, we can construct Page matrices $U(\alpha), Y(\alpha)$. In Lines 4 and 5, we obtain the submatrices  $U_{ \mathrm{p}}(\alpha), U_{ \mathrm{f}}(\alpha),\widehat{Y}_{ \mathrm{p}}(\alpha),\widehat{Y}_{ \mathrm{f}}(\alpha)$ and $\widehat{H}(\alpha)$. Due to the linearity of the system and the zero initial state, we have $\overline{H}(\alpha) = \alpha \overline{H}(1)$ and thus $\sigma_{\min}(\overline{H}(\alpha))$ is proportional to $\alpha$ if $\sigma_{\min}(\overline{H}(1))\neq 0$. Therefore, when observing that $\sigma_{\min}(\widehat{H}(1))$ increases approximately proportionally with $\alpha$, we claim that $\sigma_{\min}(\overline{H}(1))\neq 0$.
\end{remark}

\begin{algorithm}[htbp!]
\caption{Data-driven observability index identification for MISO systems}
\textbf{Input:} horizon length $L$ for the Page matrices\\
\textbf{Output:} the system order $l_\mathrm{o}$
\begin{algorithmic}[1]
\State Use \eqref{eq: design_of_inputs} to generate historical data $(u_{[1,T]},\hat{y}_{[1,T]})$ to construct Page matrices $U=\mathcal{P}_{L}(u_{[1,T]})$ and $\widehat{Y}=\mathcal{P}_{L}(\hat{y}_{[1,T]})$.
\State $k\gets 1, \textsc{TER}= 0$
	\While{$\textsc{TER} = 0$}
 \State Partition $U = [U_{ \mathrm{p}}^\top \;U_{ \mathrm{f}}^\top]^\top$, $\widehat{Y} = [\widehat{Y}_{ \mathrm{p}}^\top \;\widehat{Y}_{ \mathrm{f}}^\top]^\top$ such that $U_\mathrm{p}$ has $km$ rows (i.e., $l_{\mathrm{p}} = k$)
 \State Build $\widehat{H} = [U_{ \mathrm{p}}^\top\; \widehat{Y}_{ \mathrm{p}}^\top\; U_{ \mathrm{f}}^\top]^\top$ 
\If{$\sigma_{\min}(\widehat{H})\leq l_\mathrm{h}\delta$ }
    \State $l_\mathrm{o}\gets k-1$, TER$\gets 1$
\EndIf
	\State $k\gets k+1$
	\EndWhile
\end{algorithmic}
\label{alg:observability index identification}
\end{algorithm}

\section{Robust control for regulation of MISO systems with suboptimality guarantee}
\label{sec: formulationofGu} 
In this section, we propose a data-driven robust control method for MISO systems, where we use the prediction error bounds in Theorem \ref{thm: prediction_error_bound} to calculate $c_{\textrm{worst}}(u_\mathrm{f})$, an upperbound to the regulation cost. To justify that this upperbound is not too conservative, we study the suboptimality of the derived input sequence and compare it with the optimal input sequence. The extension to multiple-output systems is provided
in Section \ref{sec: computation+extensions}. 

With the notation
$\Delta:=\{\Delta_{Y_{\mathrm{p}}},\Delta_{Y_{\mathrm{f}}},\Delta_{y_{\mathrm{p}}}\},$ the robust regulation problem is formulated as the following bilevel
program where the inner problem calculates an upperbound for the worst-case cost and the outer problem optimizes the input such that the cost upperbound is minimized, \begin{small} 
\begin{subequations}
\label{eq:robust_optimization} 
\begin{align}
\underset{u_{\mathrm{f}}}{\text{min}}  \underset{\Delta,g,y_{\mathrm{f}}}{\text{max}}\quad & y_{\mathrm{f}}^{\top}y_{\mathrm{f}}+u_{\mathrm{f}}^{\top}u_{\mathrm{f}}\\
\text{s.t.}\qquad & \begin{bmatrix}U_{{\mathrm{p}}}\\
\widehat{Y}_{{\mathrm{p}}}\\
U_{{\mathrm{f}}}\\
\widehat{Y}_{{\mathrm{f}}}
\end{bmatrix}g+\begin{bmatrix}0\\
\Delta_{Y_{\mathrm{p}}}\\
0\\
\Delta_{{Y}_{{\mathrm{f}}}}
\end{bmatrix}g=\begin{bmatrix}u_{{\mathrm{p}}}\\
\hat{y}_{{\mathrm{p}}}\\
u_{{\mathrm{f}}}\\
y_{\mathrm{f}}
\end{bmatrix}+\begin{bmatrix}0\\
\Delta_{{y}_{{\mathrm{p}}}}\\
0\\
0
\end{bmatrix}\label{eq: suboptimality,equality}\\
& \max\left(||\Delta_{Y_{\mathrm{p}}}||_{\max},||\Delta_{Y_{\mathrm{f}}}||_{\max},||\Delta_{y_{\mathrm{p}}}||_{\infty}\right)\leq\delta\label{eq: suboptimality,noise-bound}\\
& ||g-\hat{g}(u_{\mathrm{f}})||^{2}\leq\mathcal{C}^{2}(u_{\mathrm{f}})\delta^{2}.\label{eq: suboptimality,g-bound}
\end{align}
\end{subequations}
\end{small}where $\mathcal{C}(u_\mathrm{f})$ is defined in Theorem \ref{thm: prediction_error_bound}. We use the alternating optimization method in \cite{lu2019block} to solve problem \eqref{eq:robust_optimization}(for details see Section \ref{sec: output_trajectory prediction}). We denote the solution to the outer
problem as $\check{u}_{\mathrm{f}}$. Given any input sequence $u_\mathrm{f}$, we let $\bar{c}(u_{\mathrm{f}})$  be the resulting true regulation cost and $c_{\mathrm{worst}}(u_{\mathrm{f}})$ be the optimal objective value of the inner problem of \eqref{eq:robust_optimization}.  In the following theorem, we show that $c_{\mathrm{worst}}(u_\mathrm{f})$ is indeed an upperbound to the true regulation cost.
\begin{theorem}If Assumptions \ref{ass: measurement_noises}, \ref{ass: design_of_historical_data} and
\ref{ass: small_noise} hold, for any $u_\mathrm{f}$ we have $c_{\mathrm{worst}}({u}_{\mathrm{f}})\geq \bar{c}({u}_{\mathrm{f}}).$
\label{thm: upperbound} \end{theorem} 

\begin{proof} The noise realization
$(\overline{\Delta}_{Y_{\mathrm{p}}},\overline{\Delta}_{Y_{\mathrm{f}}},\overline{\Delta}_{y_{\mathrm{p}}})$
satisfies that 
\[
\max\left(||\overline{\Delta}_{Y_{\mathrm{p}}}||_{\max},||\overline{\Delta}_{Y_{\mathrm{f}}}||_{\max},||\overline{\Delta}_{y_{\mathrm{p}}}||_{\infty}\right)\leq\delta.
\]
With the vector $\bar{g}(u_\mathrm{f})$ we can reconstruct the noiseless system output, i.e., $(\widehat{Y}_\mathrm{p}+\overline{\Delta}_{y_{\mathrm{p}}})\bar{g}(u_\mathrm{f}) = y_\mathrm{p}$, $(\widehat{Y}_\mathrm{f}+\overline{\Delta}_{y_{\mathrm{f}}})\bar{g}(u_\mathrm{f}) = \bar{y}_\mathrm{f}(u_\mathrm{f})$. Now we see that ${u}_{\mathrm{f}},\overline{\Delta}_{Y_{\mathrm{p}}},\overline{\Delta}_{Y_{\mathrm{f}}}$,
$\overline{\Delta}_{y_{\mathrm{p}}},\bar{{g}}({u}_{\mathrm{f}}),$ and $\bar{y}_{\mathrm{f}}({u}_{\mathrm{f}})$ satisfy \eqref{eq: suboptimality,equality}, \eqref{eq: suboptimality,noise-bound}
and the corresponding regulation cost is $\bar{c}(u_{\mathrm{f}})$.
Meanwhile, considering the error bounds in Theorem \ref{thm: prediction_error_bound}, the constraint \eqref{eq: suboptimality,g-bound} is also satisfied. Since $c_{\mathrm{worst}}(u_{\mathrm{f}})$ is the maximum
cost in the inner problem of \eqref{eq:robust_optimization}, we have
$c_{\mathrm{worst}}(u_{\mathrm{f}})\geq \bar{c}(u_{\mathrm{f}})$.
\end{proof} In \eqref{eq:robust_optimization}, the outer problem seeks a $u_{\mathrm{f}}$ that minimizes
this upperbound. Similar ideas that minimize a cost upperbound can be found in \cite{dean2020sample,furieri2022near}. To discuss the conservativeness of this approach, we need to bound the suboptimality
of the solution to \eqref{eq:robust_optimization} when compared to
the noiseless case. To this aim, we define the tuple $\widehat{\mathcal{Y}}:=(\widehat{Y}_{\mathrm{p}},\widehat{Y}_{\mathrm{f}},\widehat{y}_{\mathrm{p}})$, consider the control input sequences 
\begin{equation}
\label{eq: hatufbaruf}
\begin{aligned}
 \hat{u}^*_\mathrm{f}: &= \mathrm{argmin}_{u_\mathrm{f}} \|\hat{y}_\mathrm{f}(u_\mathrm{f})\|^2 + \|u_\mathrm{f}\|^2,\\
  {u}^*_\mathrm{f}: &= \mathrm{argmin}_{u_\mathrm{f}} \|\bar{y}_\mathrm{f}(u_\mathrm{f})\|^2 + \|u_\mathrm{f}\|^2,
\end{aligned}
\end{equation}
and bound the error $\|\hat{u}^*_\mathrm{f}-{u}^*_\mathrm{f}\|$ in the following lemma whose proof is given in Appendix \ref{sec: app_eta}.

\begin{lemma} \label{lmm: eta} Under Assumption \ref{ass: measurement_noises}, \ref{ass: design_of_historical_data} and \ref{ass: small_noise}, we define the following polynomials in $\delta$, 
\begin{equation}
\begin{aligned}
\mathcal{F}_1(\delta,\widehat{\mathcal{Y}}): &= 2l_{\mathrm{h}}||\widehat{H}^{\dagger}||(1+4||\widehat{Y}_{\mathrm{p}}||\cdot||\widehat{H}^{\dagger}||)\delta\\
\mathcal{F}_2(\delta,\widehat{\mathcal{Y}}): &=(2||\widehat{K}_{1}||+\mathcal{F}_{1}(\delta,\widehat{\mathcal{Y}}))\mathcal{F}_{1}(\delta,\widehat{\mathcal{Y}})\\
\mathcal{F}_3(\delta,\widehat{\mathcal{Y}}): &=\sqrt{l_{\mathrm{p}}}\left\|\widehat{K}_{2}^{\top}\widehat{K}_{1}\right\|\delta+\left(\left\|\hat{b}(0)\right\|+\sqrt{l_{\mathrm{p}}}\delta\right)\mathcal{F}_{2}(\delta,\widehat{\mathcal{Y}})\\
\mathcal{F}(\delta,\widehat{\mathcal{Y}}): &=\left\|(\widehat{K}_{2}^{\top}\widehat{K}_{2}+I)^{-1}\right\|\mathcal{F}_{3}(\delta,\widehat{\mathcal{Y}})+\left(\left\|\widehat{K}_{2}^{\top}\widehat{K}_{1}\hat{b}(0)\right\|+\mathcal{F}_{3}(\delta,\widehat{\mathcal{Y}})\right)\mathcal{F}_{2}(\delta,\widehat{\mathcal{Y}}),
\end{aligned}
\end{equation}
where $\widehat{K}_{1}:=\widehat{Y}_{\mathrm{f}}\widehat{H}^{\dagger}$,
$\widehat{K}_{2}:=\widehat{K}_{1}[0 \; 0 \; I]^{\top}$ and $\hat{b}(u_\mathrm{f})$ is defined in \eqref{eq: ls_pred_noisy}. Then, by letting $\eta(\delta,\widehat{\mathcal{Y}}) =1+\frac{||\widehat{H}^{\dagger}||\mathcal{F}(\delta,\widehat{\mathcal{Y}})}{||\hat{g}(\hat{u}^*_{\mathrm{f}})||}$, we have $||\hat{u}^*_{\mathrm{f}}-u_{\mathrm{f}}^{*}||\leq\mathcal{F}(\delta,\widehat{\mathcal{Y}})$ and
\begin{equation}
||\hat{g}({u}_{\mathrm{f}}^{*})||\leq\eta(\delta,\widehat{\mathcal{Y}})||\hat{g}(\hat{u}^*_{\mathrm{f}})||.\label{eq: multiple_inequality}
\end{equation}
\end{lemma}

Recall that $c_{\mathrm{worst}}(u_{\mathrm{f}})$ and $\bar{c}(u_{\mathrm{f}})$
are, respectively, the worst-case cost derived by the inner problem
of \eqref{eq:robust_optimization} and the resulting true regulation
cost when $u_{\mathrm{f}}$
is applied, $u_{\mathrm{f}}^{*}$ is the optimal input for the noiseless
case (see \eqref{eq: hatufbaruf}) and $\check{u}$ is the optimal solution to \eqref{eq:robust_optimization}.
In the following, we compare $\bar{c}(\check{u}_{\mathrm{f}})$
with $\bar{c}(u_{\mathrm{f}}^{*})$ to see how much suboptimality
is introduced by solving \eqref{eq:robust_optimization} and applying
$\check{u}_{\mathrm{f}}$. \begin{theorem} \label{thm: suboptimality}
Let Assumption \ref{ass: measurement_noises}, \ref{ass: design_of_historical_data} and \ref{ass: small_noise} hold and define

\begin{equation}
\begin{aligned}
\mathcal{C}_1(\delta,\widehat{\mathcal{Y}}): &= 2\sigma_{\min}^{-1}(\widehat{H})(\sqrt{l_{\mathrm{p}}}+\eta(\delta,\widehat{\mathcal{Y}})l_{\mathrm{h}}||\hat{g}(\hat{u}^*_{\mathrm{f}})||)\delta\\
\mathcal{C}_2(\delta,\widehat{\mathcal{Y}}): &=(||\widehat{Y}_{\mathrm{f}}||+l_{\mathrm{h}}\delta)\mathcal{C}_{1}(\delta,\widehat{\mathcal{Y}})+\eta(\delta,\widehat{\mathcal{Y}}) l_{\mathrm{h}}||\hat{g}(\hat{u}^*_{\mathrm{f}})||\delta\\
\mathcal{C}_3(\delta,\widehat{\mathcal{Y}}): &=8 (\mathcal{C}_2(\delta,\widehat{\mathcal{Y}}))^{2} + 4\left\|\widehat{Y}_{\mathrm{f}}\right\|\cdot||\hat{g}(\hat{u}^*_{\mathrm{f}})||\eta(\delta,\widehat{\mathcal{Y}})\mathcal{C}_2(\delta,\widehat{\mathcal{Y}}).
\end{aligned}
\end{equation}
One has that $\bar{c}(\check{u}_{\mathrm{f}})-\bar{c}(u_{\mathrm{f}}^{*})\leq\mathcal{C}_3(\delta,\widehat{\mathcal{Y}})$. Morever, the upperbound $\mathcal{C}_3(\delta,\widehat{\mathcal{Y}})$, computable by using only the noisy measurements, converges in probability to 0 as $\delta$ converges to 0.
\end{theorem}
\begin{proof} We denote $(\widetilde{\Delta}_{Y_{\mathrm{p}}},\widetilde{\Delta}_{Y_{\mathrm{f}}},\widetilde{\Delta}_{y_{\mathrm{p}}},\tilde{g}(u_{\mathrm{f}}))$ as the solution to the following optimization problem 
\begin{equation}
\begin{aligned}\underset{\Delta_{Y_{\mathrm{p}}},\Delta_{Y_{\mathrm{f}}},\Delta_{y_{\mathrm{p}}},g}{\mathrm{max}} & g^{\top}({\widehat{Y}_{\mathrm{f}}}+\Delta_{Y_{\mathrm{f}}})^{\top}({\widehat{Y}_{\mathrm{f}}}+\Delta_{Y_{\mathrm{f}}})g+u_{\mathrm{f}}^{\top}u_{\mathrm{f}}\\
\text{s.t.}\qquad & \begin{bmatrix}U_{{\mathrm{p}}}\\
\widehat{Y}_{{\mathrm{p}}}\\
U_{{\mathrm{f}}}
\end{bmatrix}g+\begin{bmatrix}0\\
\Delta_{Y_{\mathrm{p}}}\\
0
\end{bmatrix}g=\begin{bmatrix}u_{{\mathrm{p}}}\\
\hat{y}_{{\mathrm{p}}}\\
u_{{\mathrm{f}}}^{*}
\end{bmatrix}+\begin{bmatrix}0\\
\Delta_{{y}_{{\mathrm{p}}}}\\
0
\end{bmatrix}\\
& \max\left(||\Delta_{Y_{\mathrm{p}}}||_{\max},||\Delta_{Y_{\mathrm{f}}}||_{\max},||\Delta_{y_{\mathrm{p}}}||_{\infty}\right)\leq\delta\\
& ||g-\hat{g}(u_{\mathrm{f}})||^{2}\leq\mathcal{C}^{2}(u_{\mathrm{f}})\delta^{2}
\end{aligned}
\label{eq:innermaximization}
\end{equation}
and also define $\tilde{y}_{\mathrm{f}}(u_{\mathrm{f}}):=({\widehat{Y}_{\mathrm{f}}}+\widetilde{\Delta}_{y_\mathrm{f}})\tilde{g}(u_{\mathrm{f}})$.
Since $||\hat{g}(u_{\mathrm{f}}^{*})-\bar{g}(u_{\mathrm{f}}^{*})||\leq \mathcal{C}(u_{\mathrm{f}}^{*})\delta$ according to Theorem \ref{thm: prediction_error_bound}, we have
\begin{equation}
\begin{aligned}\frac{1}{2}\sigma_{\min}(\widehat{H})||\hat{g}(u_{\mathrm{f}}^{*})-\bar{g}(u_{\mathrm{f}}^{*})|| & \leq(\sqrt{l_{\mathrm{p}}}+l_{\mathrm{h}}||\hat{g}(u_{\mathrm{f}}^{*})||)\delta\\
& \leq(\sqrt{l_{\mathrm{p}}}+\eta(\delta,\widehat{\mathcal{Y}})l_{\mathrm{h}}||\hat{g}(\hat{u}^*_{\mathrm{f}})||)\delta,
\end{aligned}
\label{eq:gbound}
\end{equation}
where $\eta(\delta,\widehat{\mathcal{Y}})$ is derived in Lemma \ref{lmm: eta}. The inequality \eqref{eq:gbound} implies $||\hat{g}(u_{\mathrm{f}}^{*})-\bar{g}(u_{\mathrm{f}}^{*})||\leq\mathcal{C}_{1}(\delta,\widehat{\mathcal{Y}})$.
Similarly, we have $||\hat{g}(u_{\mathrm{f}}^{*})-\tilde{g}(u_{\mathrm{f}}^{*})||\leq\mathcal{C}_{1}(\delta,\widehat{\mathcal{Y}})$.
Consequently, we can derive \begin{small} 
\begin{equation}
\begin{aligned}||\hat{y}_{\mathrm{f}}(u_{\mathrm{f}}^{*})-\bar{y}_{\mathrm{f}}(u_{\mathrm{f}}^{*})|| & =||\widehat{Y}_{\mathrm{f}}(\hat{g}(u_{\mathrm{f}}^{*})-\bar{g}(u_{\mathrm{f}}^{*}))-\overline{\Delta}_{Y_{\mathrm{f}}}\bar{g}(u_{\mathrm{f}}^{*})||\\
& \leq||\widehat{Y}_{\mathrm{f}}||\cdot||\hat{g}(u_{\mathrm{f}}^{*})-\bar{g}(u_{\mathrm{f}}^{*})||+||\overline{\Delta}_{Y_{\mathrm{f}}}||\cdot||\bar{g}(u_{\mathrm{f}}^{*})||\\
& \leq||\widehat{Y}_{\mathrm{f}}||\cdot||\hat{g}(u_{\mathrm{f}}^{*})-\bar{g}(u_{\mathrm{f}}^{*})||+\\
& \quad\quad l_{\mathrm{h}}\delta(||\hat{g}(u_{\mathrm{f}}^{*})||+||\bar{g}(u_{\mathrm{f}}^{*})-\hat{g}(u_{\mathrm{f}}^{*})||)\\
& \leq\underbrace{(||\widehat{Y}_{\mathrm{f}}||+l_{\mathrm{h}}\delta)\mathcal{C}_{1}(\delta,\widehat{\mathcal{Y}})+\eta(\delta,\widehat{\mathcal{Y}}) l_{\mathrm{h}}||\hat{g}(\hat{u}^*_{\mathrm{f}})||\delta}_{\mathcal{C}_2(\delta,\widehat{\mathcal{Y}})},
\end{aligned}
%\label{eq:yfbound}
\end{equation}
\end{small}$||\hat{y}(u_{\mathrm{f}}^{*})-\tilde{y}(u_{\mathrm{f}}^{*})||\leq\mathcal{C}_2(\delta,\widehat{\mathcal{Y}})$
and $||\bar{y}(u_{\mathrm{f}}^{*})-\tilde{y}(u_{\mathrm{f}}^{*})||\leq2\mathcal{C}_2(\delta,\widehat{\mathcal{Y}})$.
Finally, since $\check{u}_{\mathrm{f}}$ is the solution to the outer problem
of \eqref{eq:robust_optimization} while ${u}_{\mathrm{f}}^{*}$ is feasible, $c_{\mathrm{worst}}(\check{u}_{\mathrm{f}})\leq c_{\mathrm{worst}}({u}_{\mathrm{f}}^{*})=||\tilde{y}_{\mathrm{f}}({u}_{\mathrm{f}}^{*})||^{2}+||u_{\mathrm{f}}^{*}||^{2}$, based on which we have 
\begin{equation}
\begin{aligned}\bar{c}(\check{u}_{\mathrm{f}})-\bar{c}({u}_{\mathrm{f}}^{*}) & \leq c_{\mathrm{worst}}(\check{u}_{\mathrm{f}})-\bar{c}({u}_{\mathrm{f}}^{*})\\
& \leq c_{\mathrm{worst}}({u}_{\mathrm{f}}^{*})-\bar{c}({u}_{\mathrm{f}}^{*})\\
& =\left\|\tilde{y}_{\mathrm{f}}(u_{\mathrm{f}}^{*})\right\|^{2}+||u_{\mathrm{f}}^{*}||^{2}-\bar{c}({u}_{\mathrm{f}}^{*})\\
& \leq||\bar{y}_{\mathrm{f}}(u_{\mathrm{f}}^{*})||^{2}+||u_{\mathrm{f}}^{*}||^{2}+||\tilde{y}_{\mathrm{f}}(u_{\mathrm{f}}^{*})-\bar{y}_{\mathrm{f}}(u_{\mathrm{f}}^{*})||^{2}-\bar{c}({u}_{\mathrm{f}}^{*})+\\
& \quad\text{ }2(||\hat{y}_{\mathrm{f}}(u_{\mathrm{f}}^{*})||+||\bar{y}_{\mathrm{f}}(u_{\mathrm{f}}^{*})-\hat{y}_{\mathrm{f}}(u_{\mathrm{f}}^{*})||)\cdot||\tilde{y}_{\mathrm{f}}(u_{\mathrm{f}}^{*})-\bar{y}_{\mathrm{f}}(u_{\mathrm{f}}^{*})||\\
& \leq2(\eta(\delta,\widehat{\mathcal{Y}})\left\|\widehat{Y}_{\mathrm{f}}\right\|\cdot||\hat{g}(\hat{u}^*_{\mathrm{f}})||+\mathcal{C}_2(\delta,\widehat{\mathcal{Y}}))\cdot 2\mathcal{C}_2(\delta,\widehat{\mathcal{Y}})+4(\mathcal{C}_2(\delta,\widehat{\mathcal{Y}}))^{2}\\
&= \underbrace{8 (\mathcal{C}_2(\delta,\widehat{\mathcal{Y}}))^{2} + 4\left\|\widehat{Y}_{\mathrm{f}}\right\|\cdot||\hat{g}(\hat{u}^*_{\mathrm{f}})||\eta(\delta,\widehat{\mathcal{Y}})\mathcal{C}_2(\delta,\widehat{\mathcal{Y}})}_{\mathcal{C}_3(\delta,\widehat{\mathcal{Y}})}.
\end{aligned}
\end{equation}

As $\delta$ goes to 0, the noisy measurement $\widehat{\mathcal{Y}}$ converges in probability to $\overline{\mathcal{Y}}:=({Y}_{\mathrm{p}},{Y}_{\mathrm{f}}$, $y_{\mathrm{p}})$. Therefore, we see that $\sigma_{\min}^{-1}(\widehat{H})$ converges to in probability to $\sigma_{\min}^{-1}(\overline{H})$, $\eta(\delta,\widehat{\mathcal{Y}})$ to 1 and thus $\mathcal{C}_3(\delta,\widehat{\mathcal{Y}})$ to 0.
\end{proof}
Now, we see that as the measurement noise diminishes, the regulation cost resulting from the solution to \eqref{eq:robust_optimization} decreases to the minimal value. This is different from the suboptimality bound in \cite{huang2021robust} where the achieved regulation cost is only shown to be less than twice the minimal value.
\begin{remark}
\label{rmk: averaging}
If the measurement noise sequence is i.i.d., by sampling the same historical trajectories for $N$ times to construct $\widehat{Y}_{\mathrm{p},i},\widehat{Y}_{\mathrm{f},i},\hat{y}_{\mathrm{p},i}$ for $i=1,\ldots, N$ (representing instances of ${Y}_{\mathrm{p}},{Y}_{\mathrm{f}},{y}_{\mathrm{p}}$ with independent noise realizations) and calculating the average values $\widehat{Y}^{\mathsf{avg}}_\mathrm{p} = (1/N)\sum^N_{i=1}\widehat{Y}_{\mathrm{p},i},\widehat{Y}^{\mathsf{avg}}_\mathrm{f}= (1/N)\sum^N_{i=1}\widehat{Y}_{\mathrm{f},i},\hat{y}^{\mathsf{avg}}_\mathrm{p}= (1/N)\sum^N_{i=1}\hat{y}_{\mathrm{f},i}$, we can attenuate the influence of noise. Specifically, given any $\epsilon>0$, there exists $0<\delta^{\mathsf{new}}(N,\epsilon)<\delta$ such that $\delta^{\mathsf{new}}(N,\epsilon)$ converges to 0 as $N$ goes to infinity and
$$\mathrm{max}\{||\widehat{Y}^{\mathsf{avg}}_\mathrm{p} - {Y}_\mathrm{p}||_\mathrm{max},||\widehat{Y}^{\mathsf{avg}}_\mathrm{f} - {Y}_\mathrm{f}||_\mathrm{max},||\hat{y}^{\mathsf{avg}}_\mathrm{p} - {y}_\mathrm{p}||_\mathrm{max}\}\leq \delta^{\mathsf{new}}(N,\epsilon)$$
holds with a probability of $1-\epsilon$. The averaging technique can be used to satisfy Assumption \ref{ass: small_noise}. Moreover, it allows one to conduct a sampling complexity analysis for the robust control scheme \eqref{eq:robust_optimization} (i.e., upperbounding the number of samples required to achieve a given suboptimality level), similar with the ones conducted in  
\cite{dean2020sample,xu2021finite} for model-based schemes. 
\end{remark}
\begin{remark}
If the formulation \eqref{eq:robust_optimization} is extended to solve a trajectory tracking problem with an objective function $\|y_\mathrm{f}-y_\mathrm{ref}\|_Q+\|u_\mathrm{f}\|_R$ with positive semidefinite matrices $Q\in \mathbb{R}^{p\times p}$ and $R\in \mathbb{R}^{m\times m}$, it is easy to modify Theorem \ref{thm: suboptimality} for upperbounding the suboptimality gap.
\end{remark}

\section{Extensions}
\label{sec: computation+extensions}

\subsection{Extensions to MIMO systems}
\label{sec: extensions for MIMO cases}
According to Remark \ref{rmk: mimonotgood}, for MIMO systems we cannot bound the prediction error $||\hat{y}(u_\mathrm{f}) -\bar{y}(u_\mathrm{f})||$ resulting from the scheme \eqref{eq: ls_pred_noisy} and thus cannot analyse the suboptimality of our robust control framework \eqref{eq:robust_optimization}. Here, we propose a method where Page matrices are built seperately for each output.

Suppose we have outputs $y^{1},y^{2},\ldots,y^{p}$. For each output $y^{i}$, we have a MISO sub-system with minimal realization
\begin{equation}\label{eq:subsystem} 
\begin{aligned}
x^{i}_{t+1} &=A^{i} x^{i}_{t}+B^{i} u_{t}, \\
y^{i}_{t} &=C^{i} x^{i}_{t}+D^{i} u_{t},
\end{aligned}
\end{equation}
where the matrices $A^{i},B^{i},C^{i},D^{i}$ are unknown. We build, according to \eqref{eq: design_of_inputs}, the Page matrices $U_{\mathrm{p}}^{i}$, $U_{\mathrm{f}}^{i}$, $\widehat{Y}_{\mathrm{p}}^{i}$, $\widehat{Y}_{\mathrm{f}}^{i}$, $\widehat{H}^{i}$ along with the recent vectors $u_{\mathrm{p}}^{i}$, $u_{\mathrm{f}}^{i}$, $\hat{y}_{\mathrm{p}}^{i}$ for each $i$. Similar with \eqref{eq:robust_optimization}, we can write the robust control problem as 
\begin{subequations}
\label{eq:robust_optimization_subsystems} 
\begin{align}
\underset{u_{\mathrm{f}}}{\text{min}}  \underset{\Delta,g,y_{\mathrm{f}}}{\text{max}}\quad   & ||u_{\mathrm{f}}||^{2} + \sum^{n}_{i=1}||y^{i}_{\mathrm{f}}||^{2}\\
 \text{s.t.}\quad & \forall i, \quad  \begin{bmatrix}U^{i}_{{\mathrm{p}}}\\
\widehat{Y}^{i}_{{\mathrm{p}}}\\
U^{i}_{{\mathrm{f}}}\\
\widehat{Y}^{i}_{{\mathrm{f}}}
\end{bmatrix}g^{i}+\begin{bmatrix}0\\
\Delta^{i}_{Y_{\mathrm{p}}}\\
0\\
\Delta^{i}_{{Y}_{{\mathrm{f}}}}
\end{bmatrix}g^{i}=\begin{bmatrix}u_{{\mathrm{p}}}^{i}\\
\hat{y}^{i}_{{\mathrm{p}}}\\
u^{i}_{{\mathrm{f}}}\\
y^{i}_{\mathrm{f}}
\end{bmatrix}+\begin{bmatrix}0\\
\Delta^{i}_{{y}_{{\mathrm{p}}}}\\
0\\
0
\end{bmatrix}\label{eq: suboptimality,equality,subsystem}\\
 & \max\left(||\Delta^{i}_{Y_{\mathrm{p}}}||_{\max},||\Delta^{i}_{Y_{\mathrm{f}}}||_{\max},||\Delta^{i}_{y_{\mathrm{p}}}||_{\infty}\right)\leq\delta\label{eq: suboptimality,noise-bound,subsystem}\\
 & ||g^{i}-\hat{g}^{i}(u_{\mathrm{f}})||^{2}\leq(\mathcal{C}^{i}(u_{\mathrm{f}})\delta)^{2},\label{eq: suboptimality,g-bound,subsystem}
\end{align}
\end{subequations}
where $\mathcal{C}^{i}(u_{\mathrm{f}})=2\sigma_{\min}^{-1}(\widehat{H}^{i})(\sqrt{l^{i}_{\mathrm{p}}}+l^{i}_{\mathrm{h}}||\hat{g}^{i}(u_{\mathrm{f}})||).$

To justify this formulation, we notice that, Theorem \ref{thm: upperbound} and \ref{thm: suboptimality} can be applied to every sub-system in \eqref{eq:subsystem}.
Summing up all the suboptimality bounds, we can get the suboptimality bound for the whole system. The alternating method in \cite{lu2019block} is also applicable to solve \eqref{eq:robust_optimization_subsystems}.

\subsection{Extensions to account for input and output constraints}
\label{sec: extensions for output constraints}
Suppose now there are an input constraint $u_{\mathrm{f}} \in \mathcal{U}$ and an output constraint $(\bar{y}^{1}_{\mathrm{f}}(u_{\mathrm{f}}),\ldots,\bar{y}^{p}_{\mathrm{f}}(u_{\mathrm{f}}) )\in \mathcal{Y}$. The input constraint $u_{\mathrm{f}} \in \mathcal{U}$ can be easily embedded into the outer optimization problem in \eqref{eq:robust_optimization}. To ensure satisfaction of the output constraint, we have to consider the prediction error and regard as infeasible the sequence $u_{\mathrm{f}}$ that has the slightest chance of violating the output constraint. By applying the inequalities in Theorem \ref{thm: prediction_error_bound} to every MISO subsystem, we derive the prediction error bound $\mathcal{E}^i(u_\mathrm{f},\delta)$ for the $i$-th output, i.e., $\|\hat{y}^{i}_{\mathrm{f}}(u_{\mathrm{f}}) - \bar{y}^{i}_{\mathrm{f}}(u_{\mathrm{f}}) \|\leq \mathcal{E}^i(u_\mathrm{f},\delta)$. Therefore, we can guarantee that $\bar{y}^{i}_{\mathrm{f}}(u_{\mathrm{f}})\in \mathcal{B}(\hat{y}^{i}_{\mathrm{f}}(u_{\mathrm{f}}),\mathcal{E}^{i}(u_{\mathrm{f}},\delta))$. Let $$\mathcal{Y}_\mathrm{f}(u_\mathrm{f}): = \{(y^{1}_\mathrm{f},\ldots,y^{p}_\mathrm{f})|y^{i}_\mathrm{f}\in \mathcal{B}(\hat{y}^{i}_{\mathrm{f}}(u_{\mathrm{f}}),\mathcal{E}^{i}(u_{\mathrm{f}},\delta)),\forall i\}$$ be the region where the output sequence might lie. Then, we can ensure the satisfaction of the output constraints by enforcing the input in \eqref{eq:robust_optimization_subsystems} to additionally verify that
\begin{equation}
\label{eq: robust_constraints}
\mathcal{Y}_\mathrm{f}(u_\mathrm{f})\subset \mathcal{Y}.
\end{equation} 
We call the optimization problem \eqref{eq:robust_optimization_subsystems} under \eqref{eq: robust_constraints} {Safe Data-Driven Minmax Control} (SDDMC). Enforcing \eqref{eq: robust_constraints} for general $\mathcal{Y}$ can be challenging \cite{guthrie2022closed}. However, if 
\[
\mathcal{Y} = \{(\bar{y}^{1}_\mathrm{f},\ldots,\bar{y}^{p}_\mathrm{f})|{y}^{i}_{-}\leq \bar{y}^{i}_\mathrm{f} \leq {y}^{i}_{+},\forall i\}
\] is a box constraint, the constraint \eqref{eq: robust_constraints} translates to another box constraint, which is for any $i$,
\begin{equation}
\label{eq: box constraint}
{y}^{i}_{-} + \mathcal{E}^{i}(u_{\mathrm{f}},\delta) \leq \hat{y}^{i}_\mathrm{f}(u_\mathrm{f}) \leq {y}^{i}_{+}-\mathcal{E}^{i}(u_{\mathrm{f}},\delta).
\end{equation}

\section{Numerical studies}
\label{sec: numerical studies}
In this section, through numerical experiments, we elaborate on the implementation details, verify the theoretical results and test the performance of SDDMC.
\subsection{Trajectory regulation for a SISO system}
\label{sec: output_trajectory prediction}
To begin with, we consider a SISO system with the following system matrices
\begin{equation}
\label{eq: system_setting_SISO}
A = 0.99 * \begin{bmatrix}
0.7 & 0.2 & 0\\
0.3 & 0.7 &-0.1\\
0 & -0.2 & 0.8
\end{bmatrix}, \text{ } B =\begin{bmatrix}
1 \\ 2 \\ 1.5
\end{bmatrix}, \text{ }C= \begin{bmatrix}
1 & 1 & 1
\end{bmatrix}, \text{ }D=0.
\end{equation}
We assume these matrices along with the observability index $l_\mathrm{o}=3$ are unknown to us. We aim to identify the observability index $l_\mathrm{o}$, construct the relevant Page matrices, observe the trajectory prediction error and then solve a robust control problem. 

We collect historical data by exciting the system according to \eqref{eq: design_of_inputs} with $\mathcal{Q} = \mathcal{N}(0,4)$, $L=8$ and $T = 160$. The i.i.d. measurement noise in the simulation is sampled from the uniform distribution on $[-\delta,\delta]$ with $\delta = 10^{-3}$. Algorithm \ref{alg:observability index identification} returns the correct observability index $l_\mathrm{o}=3$.  Now we can derive the matrices $U_\mathrm{p},U_\mathrm{f},Y_\mathrm{p},Y_\mathrm{f}$ according to Assumption \ref{ass: design_of_historical_data}. We use \eqref{eq: ls_pred_noisy} to predict the future output with a horizon of length 3 following a fixed recent trajectory. We notice that if $\widehat{H}$ is constructed using \eqref{eq: ls_pred_noisy} and \eqref{eq: split_of_recent} with $l_\mathrm{p} = l_\mathrm{o} = 3$, the prediction error is $1.8\times 10^{-2}$. If we set $l_\mathrm{p}=4$, the error becomes $1.4\times 10^4$. Through this comparison, we see the important role of observability index identification for the prediction scheme~\eqref{eq: ls_pred_noisy}.

We then test our robust control scheme \eqref{eq:robust_optimization} on the SISO system \eqref{eq: system_setting_SISO}. With $l_\mathrm{p} = l_\mathrm{f} = 3$ and the recent trajectory $(u_\mathrm{p},y_\mathrm{p})$ where  $u_\mathrm{p} = [-5.2254,7.2684,-22.5535]^\top$ and $y_\mathrm{p} = [-1.1242,-23.7291,13.3406]^\top$, we aim to achieve a minimum regulation cost $\bar{c}(u_\mathrm{f}) = 10u^\top_\mathrm{f} * u_\mathrm{f} + \bar{y}_\mathrm{f}(u_\mathrm{f})^\top  \bar{y}_\mathrm{f}(u_\mathrm{f})$. We use the scheme \eqref{eq:robust_optimization} to address this regulation problem. 

We elaborate on how to solve \eqref{eq:robust_optimization}. To use the alternating method, we initialize the input sequence by $u_\mathrm{f} = \hat{u}^*_\mathrm{f}$ (defined in \eqref{eq: hatufbaruf}). When the input is fixed, the inner maximization problem of  \eqref{eq:robust_optimization} is solved through \textsc{IPOPT} \cite{wachter2006implementation} using the primal-dual barrier approach, where the feasible solution $(\Delta_{Y_{\mathrm{p}}},\Delta_{Y_{\mathrm{f}}},\Delta_{y_{\mathrm{p}}},g) = (0,0,\widehat{Y}_\mathrm{p}\hat{g}(u_\mathrm{f})-\hat{y}_\mathrm{p},\hat{g}(u_\mathrm{f}))$ is set as the initial point. After deriving $(\widetilde{\Delta}_{Y_{\mathrm{p}}},\widetilde{\Delta}_{Y_{\mathrm{f}}},\widetilde{\Delta}_{y_{\mathrm{p}}},\tilde{g})$, the solution to the inner problem, the outer problem of \eqref{eq:robust_optimization} is simply 
$$\min_{u_\mathrm{f}} \widetilde{G}(\widetilde{\Delta}_{Y_{\mathrm{p}}},\widetilde{\Delta}_{Y_{\mathrm{f}}},\widetilde{\Delta}_{y_{\mathrm{p}}},\tilde{g},u_{\mathrm{f}}):= \left\|(\widetilde{Y}_\mathrm{f} + \widetilde{\Delta}_{Y_{\mathrm{f}}})\begin{bmatrix} U_\mathrm{p}\\ \widehat{Y}_\mathrm{p}+ \widetilde{\Delta}_{Y_{\mathrm{p}}}\\{U}_\mathrm{f}\end{bmatrix}^\dagger\begin{bmatrix} u_\mathrm{p}\\ y_\mathrm{p}+\widetilde{\Delta}_{y_{\mathrm{p}}}\\u_\mathrm{f}\end{bmatrix}\right\|+\|u_\mathrm{f}\|^2,$$
where the objective is a quadratic function of $u_{\mathrm{f}}$ and the explicit expression of the solution can be readily computed. With the new input sequence at hand, we can start the next iteration. We terminate this algorithm when the difference of input sequences derived in two iterations has a 2-norm less than $10^{-4}$. Currently, we do not have any results on the convergence for this iterative scheme. Empirically, we observe in solving the regulation problem defined above that at most three iterations are needed before the termination even when $\delta$ increases to 1. 

\begin{figure}
\centering
\includegraphics[width = 10cm]{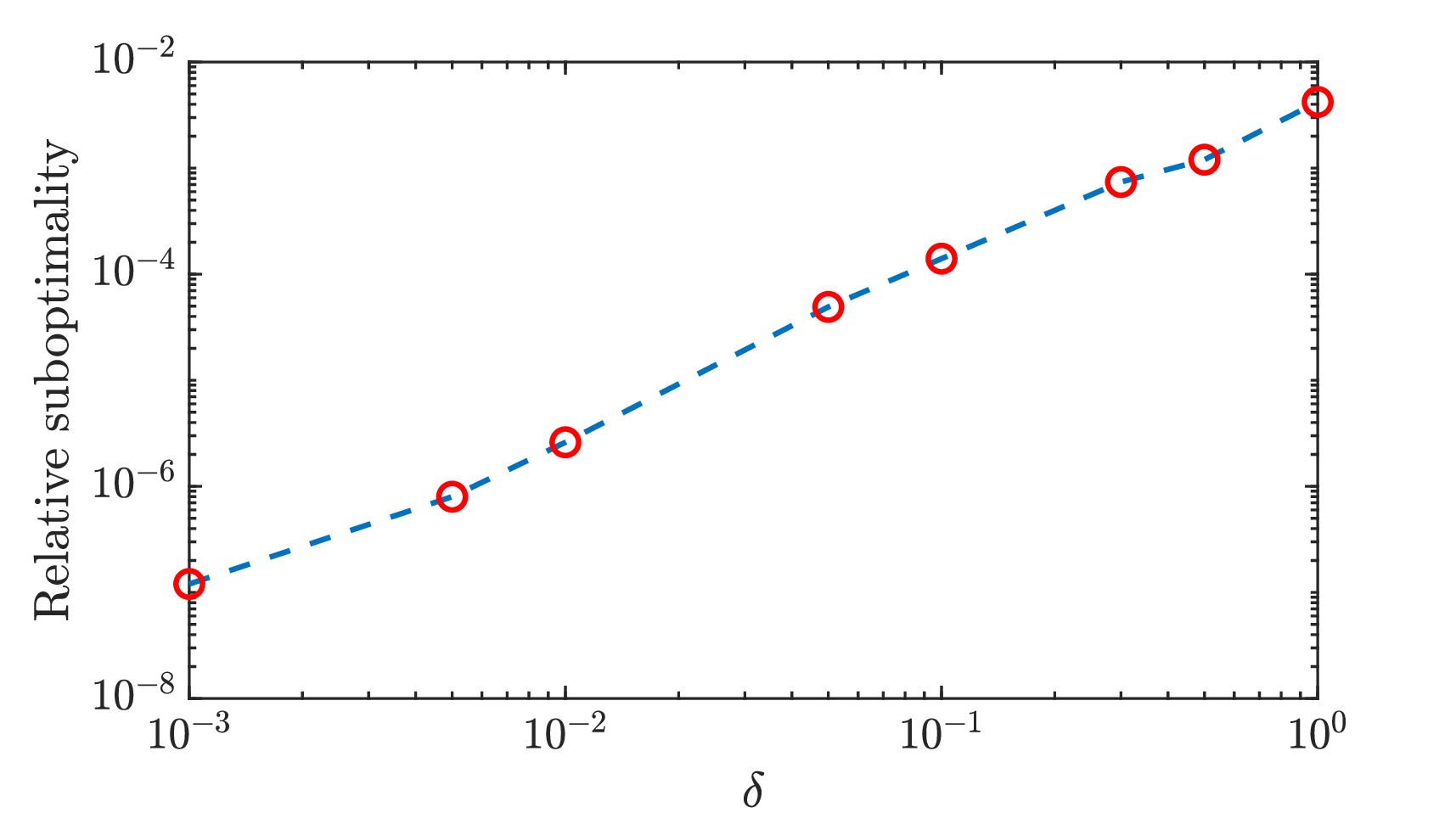}
\caption{Mean relative suboptimality for different $\delta$}
\label{fig:relative_sub}
\end{figure}

Then, we solve 50 instances of the regulation problem with independent noise realizations. To evaluate the performance of the control scheme \eqref{eq:robust_optimization}, we calculate for each instance the relative suboptimality, defined as $\frac{\bar{c}(\check{u}_{\mathrm{f}})-\bar{c}({u}^*_{\mathrm{f}})}{\bar{c}(\check{u}_{\mathrm{f}})}$. In Fig. \ref{fig:relative_sub}, we show the mean relative suboptimality for different $\delta$. We see that the input sequence derived based on \eqref{eq:robust_optimization} and the alternating method achieves a suboptimality gap that decreases to 0 as the noise diminishes, which coincides with Theorem \ref{thm: suboptimality}. Although we find that with $\delta>10^{-2}$ Assumption \ref{ass: small_noise} is not satisfied and thus Theorem \ref{thm: suboptimality} is not valid, the relative suboptimality resulting from the robust control scheme \eqref{eq: robust_constraints} is empirically small even when $\delta = 1$. When the measurement noise is i.i.d., we can weaken Assumption \ref{ass: small_noise} such that Theorem \ref{thm: full-rankness} and Theorem \ref{thm: suboptimality} are valid for larger noise level for a high probability. However, to achieve this, we need to use random matrix analysis to tighten \eqref{eq: upperbound_of_E}, which is out of the scope of this paper. 

\subsection{Application to a MIMO system with constraints: room temperature control}
We apply Safe Data-Driven Minmax Control (SDDMC) proposed in Section \ref{sec: extensions for output constraints} to room temperature control. We consider a small building model taken from \cite{oldewurtel2008tractable}, where the temperature dynamic is normalized and linearized at the equilibrium temperature $T = 15^\circ \mathrm{C}$ (which coincides with the outdoor temperature), with a sampling time of 0.5 hour. The model is described by \eqref{eq:LTI_system_ss}, with 
$$A=\left[\begin{array}{lll}0.8511 & 0.0541 & 0.0707 \\ 0.1293 & 0.8635 & 0.0055 \\ 0.0989 & 0.0032 & 0.7541\end{array}\right], \quad B=\left[\begin{array}{l}0.07 \\ 0.006 \\ 0.004\end{array}\right], \quad C=\left[\begin{array}{lll}1 & 0 & 0 \\ 0 & 1 & 0 \end{array}\right], D=0.$$
The three states represent the temperature at different spots while only two of them can be measured as indicated in the matrix $C$. Notice that, due to normalization, a state with a value $x$ translates to $(x+15)^\circ \mathrm{C}$ at the corresponding spot. We assume that, at time $t=0$, the system is at equilibrium. Therefore, the indoor and outdoor temperature is $15^\circ \mathrm{C}$. To increase the indoor temperature to $25^\circ \mathrm{C}$ and balance between user comfort and energy consumption, we let the cost function be  $\bar{c}(u_\mathrm{f}) = 10u^\top_\mathrm{f} * u_\mathrm{f} + (\bar{y}_\mathrm{f}(u_\mathrm{f})-10*[1\;1]^\top)^\top (\bar{y}_\mathrm{f}(u_\mathrm{f})-10*[1\;1]^\top)$ covering the horizon $0\leq t \leq 4$, while enforcing the the first output to be no less than $5$ for $t\geq 1$, i.e., the first spot to be at least $20^\circ\mathrm{C}$ in half an hour after the experiment starts.
\begin{figure}[htbp!]
    \centering
    \includegraphics[width=0.6\linewidth]{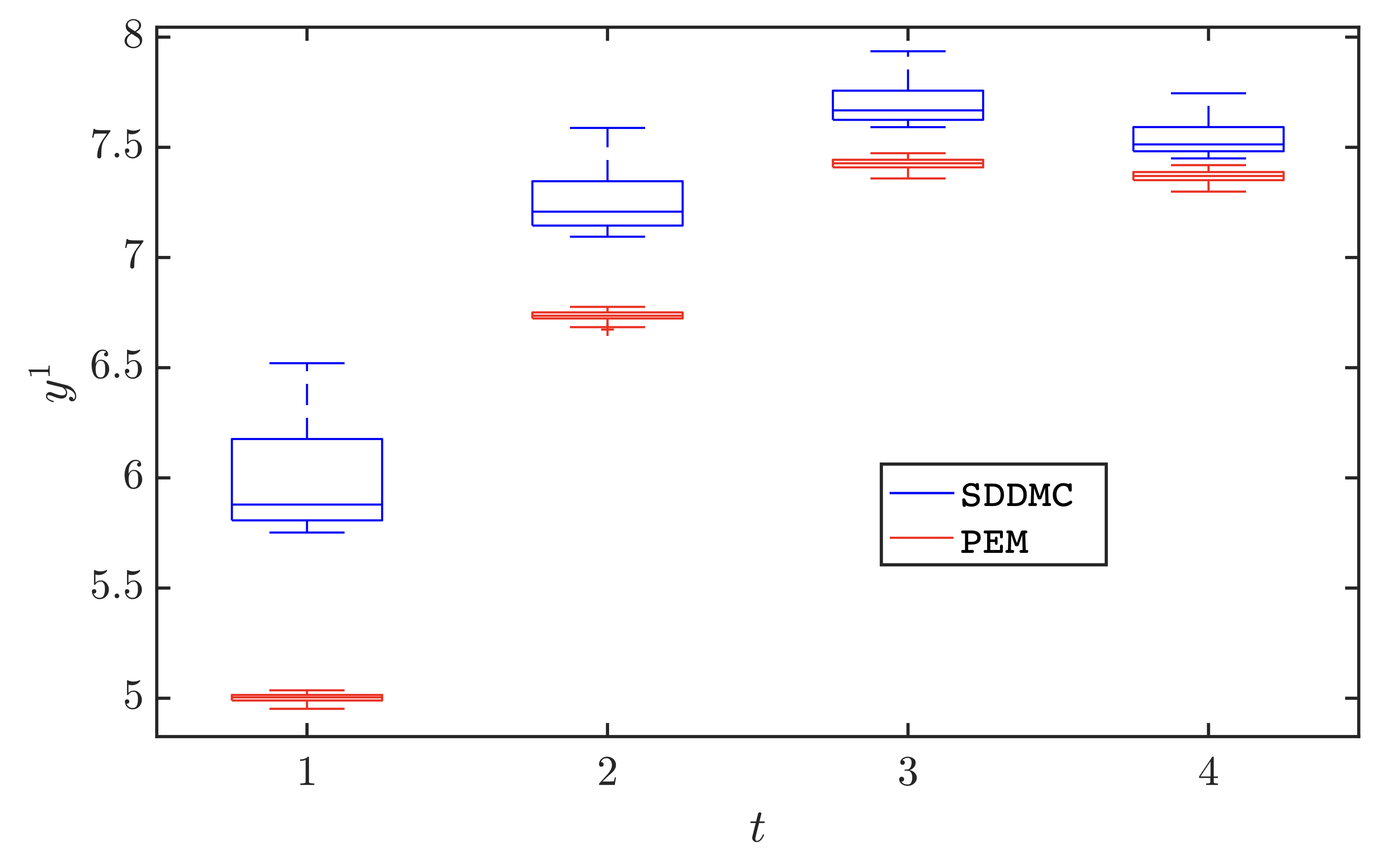}
    \caption{Trajectory of the first output: SDDMC v.s. PEM ($\delta=0.01$)}
    \label{fig:trajectory}
\end{figure}

\begin{figure}[htbp!]
    \centering
    \includegraphics[width=0.6\linewidth]{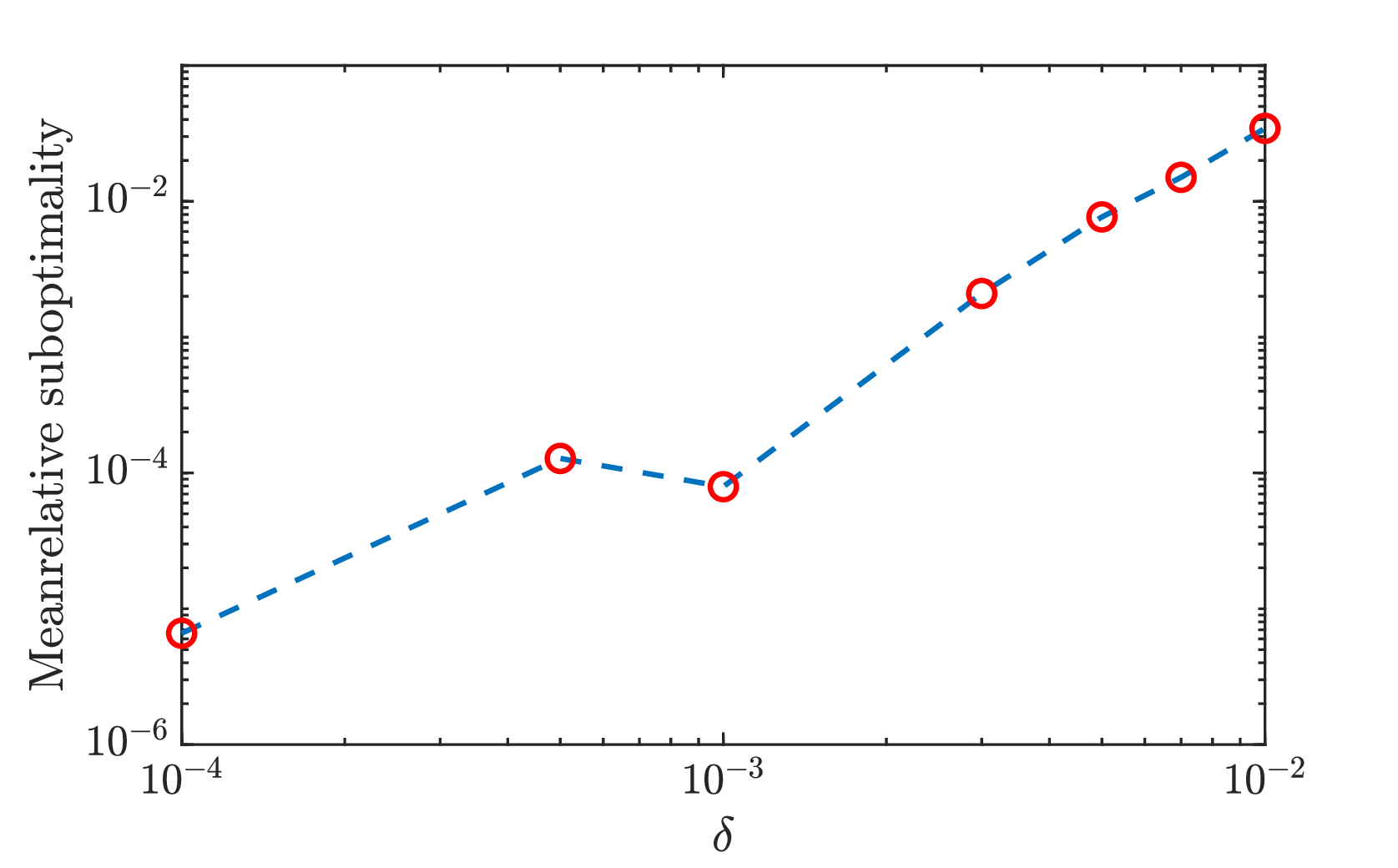}
    \caption{Room temperature control: mean relative suboptimality for different $\delta$}
    \label{fig:suboptimality_room}
\end{figure}

We run SDDMC in 50 experiments with independent noise realizations ($\delta = 0.01$). We show in Fig. \ref{fig:trajectory} the trajectory of the first output. We also compare with PEM in \cite{huang2019data} where the prediction given by \eqref{eq: ls_pred_noisy} is regarded as the true output. From Fig. \ref{fig:trajectory}, we observe that the lower bound $20^\circ \mathrm{C}$ for the first output is always satisfied if SDDMC is applied while PEM generates several trajectories with constraint violations. Meanwhile, we also plot the mean relative suboptimality in Fig.~\ref{fig:suboptimality_room} for different $\delta$. The suboptimality decreases to 0 when $\delta$ diminishes. A source of suboptimality, when $\delta$ gets larger, is the conservative estimate of the prediction error bound. This is also the main reason why in Fig. \ref{fig:trajectory} the output trajectory is substantially far away from the constraint boundaries. If one can shrink the uncertainty set, the output trajectory can get closer to the boundary, therefore leading to better optimality. 

\section{Conclusion}
\label{sec: conclusion}
In this paper, we proposed a method to construct the Page matrices for linear system trajectory prediction such that the error due to the measurement noise can be bounded using solely collected data. Based on this error bound, we designed the minimax robust control scheme such that the suboptimality gap is bounded. This scheme is extended to solve regulation problems for MIMO systems with input/output constraints. The experiments illustrated that with our method of Page matrix construction, we can achieve small prediction error and, by using the proposed robust control method, the suboptimality gap can be small for unconstrained optimal control problems. For constrained problems, we saw that the constraints were respected for different noise realizations. Our future studies will be focused on the design of historical input sequences for the derivation of a less conservative prediction error upperbound and on the development of robust control schemes with a reduced suboptimality gap.
\bibliographystyle{ieeetr}
\bibliography{references}

\appendix
\subsection{A preliminary lemma for the proof of Theorem 
\ref{thm: full-rankness}}
\label{sec: lemma_full_rankness}

\begin{lemma}
\label{lmm: full_rankness}
If the system \eqref{eq:LTI_system_ss} is SISO and Assumption \ref{ass: design_of_historical_data} holds, we let $$\mathcal{W}: = \begin{bmatrix}  
u_{L+1}  & \cdots & u_{2L}\\
u_{3L+1}  & \cdots & u_{4L}\\
\vdots & \vdots & \vdots\\
u_{2Ll_{L}-L+1} & \cdots & u_{2Ll_{L}}
\end{bmatrix},\;\;\mathcal{X}: =\begin{bmatrix}  
Cx_{L+1}  & \cdots & CA^{l_\mathrm{p}-1}x_{L+1}\\
Cx_{3L+1} & \cdots & CA^{l_\mathrm{p}-1}x_{3L+1} \\
\vdots & \vdots & \vdots\\
Cx_{2Ll_{L}-L+1}  & \cdots & CA^{l_\mathrm{p}-1}x_{2Ll_{L}-L+1} 
\end{bmatrix},$$
\begin{equation}
\label{eq: V_def}
\text{and}\quad \mathcal{V}:=\begin{bmatrix}\mathcal{W} &  \mathcal{X} \end{bmatrix}.     
\end{equation}
Then $$\mathbb{P}(\{\mathcal{V}\text{ is full rank}\}) = 1.$$
\end{lemma}

\begin{proof}
We first notice that $l_{\mathrm{p}} = l_{\mathrm{o}}\geq n_{x}$, otherwise,
for a single-output system, $\mathcal{O}(l_{\mathrm{p}})$ cannot
be column full rank. By the Cayley-Hamilton theorem, we know $\mathcal{O}(l_{\mathrm{p}})$
has the same row rank as $\mathcal{O}(n_{x})$. Since $\mathcal{O}(l_{\mathrm{p}})$
is full column rank, $\mathcal{O}(n_{x})$ is also. According to Definition
\ref{def: observability index}, $l = l_{\mathrm{p}}$ is the smallest such that $\mathcal{O}(n_{x})$ is full row rank and thus $l_{\mathrm{p}}\leq n_{x}$. Considering these facts, we conclude that $l_{\mathrm{p}}=n_{x}$. Now we proceed by proving
that the event $F_{\mathcal{X}}=\{\mathcal{X}_{1:l_{\mathrm{p}},\cdot}\text{ is full row rank}\}$
takes place almost surely.

In order to show $\mathbb{P}(F_{\mathcal{X}}) = 1$, we use an induction argument. We first look into the first row of $\mathcal{X}$. According to the system dynamic  \eqref{eq:LTI_system_ss}, we have 
\begin{equation}
x_{L+1} =\begin{bmatrix}B & AB & \cdots & A^{L-1}B\end{bmatrix}\begin{bmatrix}u_{L} & u_{L-1} & \cdots & u_{1}\end{bmatrix}^{\top}\;,\;\;\mathcal{X}_{1,\cdot}=x_{L+1}^{\top}\mathcal{O}(l_{\mathrm{p}})^{\top}.
\label{eq: non-trivial first row}
\end{equation}
Since $x_{L+1}=0$ defines a subspace in $\mathbb{R}^{mL}$ where $u_{[1,L]}$ resides, by using \eqref{eq: zero_mass} we have $\mathbb{P}\{\mathcal{X}_{1,\cdot}=0\}=\mathbb{P}\{x_{L+1}=0\}=0$.
Then, in the following, we show that, for any $l<l_{\mathrm{p}}$, 
\begin{equation}
\label{eq: induction for as}
    \mathbb{P}\{\mathcal{X}_{1:(l+1),\cdot}\text{ has row full rank}|\mathcal{X}_{1:l,\cdot}\text{ has row full rank}\}=1.
\end{equation}
To this end, we let $\Theta$ be the event where the inputs are fixed from $t=1$ to $t=(2l-1)L$ and suppose $\Theta\in\{\mathcal{X}_{1:l,\cdot}\text{ has row full rank}\}$, then there exists a vector $x_{\mathrm{n}}\in\mathbb{R}^{l_{\mathrm{p}}}$
and $x_{\mathrm{n}}\neq0$ such that $\mathcal{X}_{i,\cdot}x_{\mathrm{n}}=0$
holds for any $i\leq l$. If $\mathcal{X}_{1:(l+1),\cdot}$ does not
have full rank, we have 
\begin{equation}
\mathcal{X}_{l+1,\cdot}x_{\mathrm{n}}=0.\label{eq: orth_equation}
\end{equation}
By noticing $\mathcal{X}_{l+1,\cdot}=x_{(2l+1)L+1}^{\top}\mathcal{O}(l_{\mathrm{p}})^{\top}$, we can rewrite \eqref{eq: orth_equation} into
\begin{equation}
x_{\mathrm{n}}^{\top}\mathcal{O}(l_{\mathrm{p}})x_{(2l-1)L+1}+x_{\mathrm{n}}^{\top}\mathcal{O}(l_{\mathrm{p}})\begin{bmatrix}B & \cdots & A^{2L-1}B\end{bmatrix}\begin{bmatrix}u_{(2l+1)L} & \cdots & u_{(2l-1)L+1}\end{bmatrix}^{\top}=0.\label{eq: orth_equation2}
\end{equation}
By noticing that $\mathcal{O}(l_{\mathrm{p}})$ is full rank and $x_{\mathrm{n}}^{\top}\mathcal{O}(l_{\mathrm{p}})\neq0$,
we see $x_{\mathrm{n}}^{\top}\mathcal{O}(l_{\mathrm{p}})\begin{bmatrix}B & \cdots & A^{2L-1}B\end{bmatrix}\neq0$
since $\begin{bmatrix}B & \cdots & A^{2L-1}B\end{bmatrix}$ is full
row rank. Then, the realizations of $u_{(2l-1)L+1},\ldots,u_{(2l+1)L}$
avoids \eqref{eq: orth_equation2} almost surely due to \eqref{eq: zero_mass}. Therefore, the claim \eqref{eq: induction for as} holds because
\[
\mathbb{P}\{\mathcal{X}_{1:(l+1),\cdot}\text{ has row full rank}|\Theta\}=1.
\]Due to \eqref{eq: induction for as}, if $l<l_\mathrm{p}$ and $\mathbb{P}\{\mathcal{X}_{1:l,\cdot}\text{ has row full rank}\}=1$, then $\mathbb{P}\{\mathcal{X}_{1:{(l+1)},\cdot}\text{ has row full rank}\}=1.$ By induction, we know $\mathbb{P}\{\mathcal{X}_{1:{l_\mathrm{p}},\cdot}\text{ has row full rank}\}=1$.

Now, almost surely, the first $l_{\mathrm{p}}$ rows of $\mathcal{X}$
are linear independent while the last $L$ rows of $\mathcal{W}$ are linear independent. Since the inputs in the $(l_\mathrm{p}+i)$-th row of $\mathcal{W}$ are independent of the elements in the $(l_\mathrm{p}+i)$-th row of $\mathcal{X}$ for $1 \leq i \leq L$, we can again use the induction technique to show that the first $ (l_\mathrm{p}+i)$ rows of $ \mathcal{V}$ are linearly independent almost surely for $1 \leq i \leq L$. Thus, we have that $\mathbb{P}(\{\mathcal{V}\text{ is full rank}\})=1$.
\end{proof}

\subsection{Proof of Theorem \ref{thm: full-rankness}}
\label{sec: proof_full_rank}
In this proof, we \textbf{only} considers single-input cases for simplicity since it can be easily adapted to multiple-input cases. The main idea is to exploit the relationship between $\overline{H}$ and $\mathcal{V}$ in Lemma \ref{lmm: full_rankness}.

We let $l_L = L+l_p$ and denote $\overline{H}^\circ$ as the submatrix consisting of the $2,4,\ldots,2l_L$-th columns of $\overline{H}$. We notice that for any $i,j,$
\begin{equation}
\label{eq: y=xu}
y_{L(2i-1)+j} = CA^{j-1}x_{L(2i-1)+1} +\sum^{j-1}_{k=1}CA^kB u_{L(2i-1)+k}.   
\end{equation}  
 By substituting \eqref{eq: y=xu} into $\overline{H}^\circ$ and using elementary row transformation to simplify $\overline{H}^\circ$, we see that $\overline{H}^\circ$ is full row rank if $\mathcal{V}$, defined in \eqref{eq: V_def}, has full rank. From Lemma \ref{lmm: full_rankness}, we have that $\mathbb{P}(\{ \overline{H} \text{ is full row rank}\}) \geq \mathbb{P}(\{\overline{H}^\circ\text{ is full row rank}\})\geq \mathbb{P}(\{\mathcal{V}\text{ is full rank}\}) = 1$.

\subsection{Proof of Proposition \ref{prop: not-full-rank}}
\label{app: proof_of_not_full_row_rank}
The proof for the case $l_\mathrm{p}<l_\mathrm{o}$ is an easy modification of that for Theorem \ref{thm: full-rankness} and thus omitted. If $l_\mathrm{p}>l_\mathrm{o}$, we notice that there exist matrices $\mathcal{A}_i\in \mathbb{R}^{p\times p}$, $\mathcal{F}_i\in \mathbb{R}^{p\times m}$ for $1 \leq i \leq l_\mathrm{o}$ such that the $l_\mathrm{o}$th-order ARX model 
\begin{equation}
y_t = \sum^{l_\mathrm{o}}_{i=1}\left( \mathcal{A}_i y_{t-i} + \mathcal{F}_i u_{t-i}\right)
\end{equation}
describes exactly the system \eqref{eq:LTI_system_ss}. Therefore, we almost surely have that the $(l_\mathrm{o}p+1)$-th row of ${Y}_\mathrm{p}$ is a linear combination of the first $l_\mathrm{o}p$ rows of $U_\mathrm{p}$ and the first $l_\mathrm{o}p$ rows of ${Y}_\mathrm{p}$, which concludes the proof.

\subsection{A supporting lemma for Proposition \ref{prop: l_o_identification}}
\label{sec: supp_l_o_id}
\begin{lemma}
\label{lmm: supp_l_o_id}
Given a full-row-rank matrix $A\in \mathbb{R}^{m\times n}$ with $m<n$ and a row vector $z\in \mathbb{R}^{m\times n}$, we have that 
$$\sigma_{\min}([A^\top z^\top]^\top)\leq \sigma_{\min}(A).$$
\end{lemma}
\begin{proof}
According to the definition of singular values, we have
\begin{equation}
\begin{aligned}
\sigma_{\min}([A^\top z^\top]^\top)&=\min_{v\in\mathbb{R}^{m+1},\|v\|=1} \|[A^\top z^\top] v\|\\
&\leq \min_{v\in\mathbb{R}^{m+1},\|v\|=1, v_{m+1} = 0} \|[A^\top z^\top] v\|\\
&\leq \min_{w\in\mathbb{R}^{m},\|w\|=1} \|A^\top  v\| = \sigma_{\min}(A)
\end{aligned}
\end{equation}
\end{proof}

\subsection{Proof of Lemma \ref{lmm: eta}}
\label{sec: app_eta}

%\begin{proof}
    By substituting \eqref{eq: ls_pred_noisy} to \eqref{eq: hatufbaruf}, we have
\begin{equation}
\begin{aligned}\hat{u}^*_\mathrm{f} = \underset{u_{\mathrm{f}}}{\text{argmin}} & \text{ }(\widehat{K}_{1}\hat{b}(0)+\widehat{K}_{2}u_{\mathrm{f}})^{\top}(\widehat{K}_{1}\hat{b}(0)+\widehat{K}_{2}u_{\mathrm{f}})+u_{\mathrm{f}}^{\top}u_{\mathrm{f}}.\end{aligned}
\label{eq:CE control_reformulation}
\end{equation}
The solution is $\hat{u}^*_{\mathrm{f}}=-(\widehat{K}_{2}^{\top}\widehat{K}_{2}+I)^{-1}\widehat{K}_{2}^{\top}\widehat{K}_{1}\hat{b}(0)$.
With $\overline{K}_{1}$, $\overline{K}_{2}$ being the noiseless counterparts of $\widehat{K}_{1}$, $\widehat{K}_{2}$ respectively,
we have $u_{\mathrm{f}}^{*}=-(\overline{K}_{2}^{\top}\overline{K}_{2}+I)^{-1}\overline{K}_{2}^{\top}\overline{K}_{1}\bar{b}(0).$
In the following, we aim to bound $||\hat{u}^*_{\mathrm{f}}-u_{\mathrm{f}}^{*}||$.

Firstly, we analyse the influence of noise on $\widehat{K}_{1}$ and $\widehat{K}_{2}$.
Due to Assumption \ref{ass: small_noise},
\[
\sigma_{\min}(\overline{H})\geq\sigma_{\min}(\widehat{H})-||E||\geq\frac{1}{2}\sigma_{\min}(\widehat{H}).
\]
Thus, $||\overline{H}^{\dagger}||=\sigma_{\min}^{-1}(\overline{H})\leq2\sigma_{\min}^{-1}(\widehat{H})$.
Considering the following conclusion in perturbation
analysis \cite{stewart1977perturbation}, 
\begin{equation}
\begin{aligned}||\overline{H}^{\dagger}-\widehat{H}^{\dagger}|| & \leq2\max\{||\overline{H}^{\dagger}||^{2},||\widehat{H}^{\dagger}||^{2}\}||E||\end{aligned}
\label{eq: pseudo_inverse_perturbation}
\end{equation}
we have $||\overline{H}^{\dagger}-\widehat{H}^{\dagger}||\leq8l_{\mathrm{h}}||\widehat{H}^{\dagger}||^{2}\delta$
and 
\begin{equation}
\begin{aligned}||\widehat{K}_{1}-\overline{K}_{1}|| & \leq||\widehat{Y}_{\mathrm{p}}||\cdot||\widehat{H}^{\dagger}-{H}^{\dagger}||+||\widehat{Y}_{\mathrm{p}}-{Y}_{\mathrm{p}}||\cdot||\overline{H}^{\dagger}||\\
& \leq8l_{\mathrm{h}}||\widehat{Y}_{\mathrm{p}}||\cdot||\widehat{H}^{\dagger}||^{2}\delta+2l_{\mathrm{h}}||\widehat{H}^{\dagger}||\delta\\
& =\underbrace{2l_{\mathrm{h}}||\widehat{H}^{\dagger}||(1+4||\widehat{Y}_{\mathrm{p}}||\cdot||\widehat{H}^{\dagger}||)\delta}_{\mathcal{F}_{1}(\delta,\widehat{\mathcal{Y}})}.
\end{aligned}
\label{eq: kperturbation}
\end{equation}
Therefore, $||\widehat{K}_{2}-\overline{K}_{2}||\leq{\mathcal{F}_{1}}\delta.$ Similarly
we have 
\[
\max_{i=1,2}\left\{||\widehat{K}_{2}^{\top}\widehat{K}_{i}-\overline{K}_{2}^{\top}\overline{K}_{i}||\right\}\leq\underbrace{(2||\widehat{K}_{1}||+\mathcal{F}_{1}(\delta,\widehat{\mathcal{Y}}))\mathcal{F}_{1}(\delta,\widehat{\mathcal{Y}})}_{\mathcal{F}_{2}(\delta,\widehat{\mathcal{Y}})}.
\]

Secondly, again following \eqref{eq: pseudo_inverse_perturbation}, we have
\begin{equation}
\label{eq: second_conclusion_pseudo}
\begin{aligned}(\widehat{K}_{2}^{\top}\widehat{K}_{2}+I)^{-1}-(\overline{K}_{2}^{\top}\overline{K}_{2}+I)^{-1}\leq2||\widehat{K}_{2}^{\top}\widehat{K}_{2}-\overline{K}_{2}^{\top}\overline{K}_{2}||.\end{aligned}
\end{equation}
Through the same technique for deducing \eqref{eq: kperturbation} based on the upperbound of $\|\overline{H}^\dagger-\widehat{H}^\dagger\|$, we utilize \eqref{eq: second_conclusion_pseudo} to conclude that
\[
\begin{aligned}
\left\|\widehat{K}_{2}^{\top}\widehat{K}_{1}\hat{b}(0)-\overline{K}_{2}^{\top}\overline{K}_{1}\bar{b}(0)\right\|\leq&\underbrace{\sqrt{l_{\mathrm{p}}}\left\|\widehat{K}_{2}^{\top}\widehat{K}_{1}\right\|\delta+\left(\left\|\hat{b}(0)\right\|+\sqrt{l_{\mathrm{p}}}\delta\right)\mathcal{F}_{2}(\delta,\widehat{\mathcal{Y}})}_{\mathcal{F}_{3}(\delta,\widehat{\mathcal{Y}})},\\
\left\|\hat{u}^*_{\mathrm{f}}-u_{\mathrm{f}}^{*}\right\|\leq&\underbrace{\left\|(\widehat{K}_{2}^{\top}\widehat{K}_{2}+I)^{-1}\right\|\mathcal{F}_{3}(\delta,\widehat{\mathcal{Y}})+\left(\left\|\widehat{K}_{2}^{\top}\widehat{K}_{1}\hat{b}(0)\right\|+\mathcal{F}_{3}(\delta,\widehat{\mathcal{Y}})\right)\mathcal{F}_{2}(\delta,\widehat{\mathcal{Y}})}_{\mathcal{F}(\delta,\widehat{\mathcal{Y}})}
\end{aligned}
\]
and $\mathcal{F}(\delta,\widehat{\mathcal{Y}})$ converges to 0 as $\delta$ goes to 0.

Finally, we have $||\hat{g}(u_{\mathrm{f}}^{*})-\hat{g}(\hat{u}^*_{\mathrm{f}})||\leq||\widehat{H}^{\dagger}||\cdot||\hat{u}^*_{\mathrm{f}}-u_{\mathrm{f}}^{*}||$
and then 
\begin{equation}
\begin{aligned}||\hat{g}(u_{\mathrm{f}}^{*})||  \leq||\hat{g}(\hat{u}^*_{\mathrm{f}})||+||\widehat{H}^{\dagger}||\mathcal{F}(\delta,\widehat{\mathcal{Y}}) =||\hat{g}(\hat{u}^*_{\mathrm{f}})||(1+\frac{||\widehat{H}^{\dagger}||\mathcal{F}(\delta,\widehat{\mathcal{Y}})}{||\hat{g}(\hat{u}^*_{\mathrm{f}})||}).\end{aligned}
\end{equation} 
%\end{proof}
\end{document}